\documentclass[a4paper,10pt]{article}

\usepackage[a4paper,left=3cm,right=3cm,top=2.5cm,bottom=2.5cm]{geometry}
\usepackage{hyperref}
\usepackage{subfig}
\usepackage{pgfplots}
\usepackage{pgfplotstable}
\usepackage{mathtools}
\usepackage{multicol}
\usepackage{comment}
\usepackage{booktabs}
\pgfplotsset{compat=1.5}
\usepackage{amssymb}
\usepackage{url}
\usepackage{bm}

\usepackage{pdfsync}
\usepackage{float}
\usepackage{tabularx}
\usepackage{enumerate}
\usepackage{array}
\usepackage{xspace}
\usepackage{siunitx}
\usepackage{authblk}

\usepackage[draft,inline,marginclue]{fixme}
\usepackage{mathrsfs}
\usepackage{color, colortbl}
\usepackage{multirow}

\usepackage{amsthm}
\usepackage{amsmath}
\usepackage{stmaryrd}
\usepackage[ruled,vlined,linesnumbered]{algorithm2e}
\usepackage{graphicx}
\usepackage{booktabs}
\usepackage{cleveref}

\FXRegisterAuthor{mt}{amt}{MT} \FXRegisterAuthor{fr}{afr}{FR}
\FXRegisterAuthor{al}{aal}{AL}

\theoremstyle{remark}
\newtheorem{remark}{Remark}

\newcommand{\mupar}{\ensuremath{\boldsymbol{\mu}}}

\newcommand{\x}{\ensuremath{\mathbf{x}}}

\newcommand{\X}{\ensuremath{\mathbf{X}}}

\newcommand{\R}{\ensuremath{\mathbb{R}}}

\DeclareMathOperator*{\argmin}{\arg\!\min}
\newcolumntype{C}[1]{>{\centering\arraybackslash}m{#1}}

\definecolor{Gray}{gray}{0.9}

\newcommand{\RA}[1]{{\color{black}#1}}
\newcommand{\RB}[1]{{\color{black}#1}}

\begin{document}

\title{Multi-fidelity data fusion through parameter space reduction with
  applications to automotive engineering}
\author[a]{Francesco~Romor\footnote{francesco.romor@sissa.it}}
\author[a]{Marco~Tezzele\footnote{marco.tezzele@sissa.it}}
\author[b]{Markus~Mrosek\footnote{markus.mrosek@volkswagen.de}}
\author[b]{Carsten~Othmer\footnote{carsten.othmer@volkswagen.de}}
\author[a]{Gianluigi~Rozza\footnote{gianluigi.rozza@sissa.it}}

\affil[a]{Mathematics Area, mathLab, SISSA, via Bonomea 265, I-34136 Trieste,
  Italy}
\affil[b]{Volkswagen AG, Innovation Center Europe, 38436 Wolfsburg, Germany}

\maketitle

\begin{abstract}
  Multi-fidelity models are of great importance due to their
  capability of fusing information coming from different numerical simulations, surrogates, and sensors. We focus on the
  approximation of high-dimensional scalar functions with low
  intrinsic dimensionality. By introducing a low dimensional bias we can fight the curse
  of dimensionality affecting these quantities of interest,
  especially for many-query applications. We seek a
  gradient-based reduction of the parameter space through linear
  active subspaces or a nonlinear transformation of the input space. Then we build a
  low-fidelity response surface based on such reduction, thus enabling
  nonlinear autoregressive multi-fidelity Gaussian process regression without the need of
  running new simulations with simplified physical models. This has a
  great potential in the data scarcity regime affecting many
  engineering applications. In this
  work we present a new multi-fidelity approach that involves active subspaces and the nonlinear level-set learning method, starting from
  the preliminary analysis previously conducted~\cite{romor2020pamm}. The proposed
  framework is tested on two high-dimensional benchmark functions, and
  on a more complex car aerodynamics problem. We show how a
  low intrinsic dimensionality bias can increase the accuracy of
  Gaussian process response surfaces.
\end{abstract}


\section{Introduction}
\label{sec:intro}
The curse of dimensionality affects the realization of reliable models
for high-dimensional functions approximation. This problem is
particularly evident in the data scarcity regime which characterizes
many industrial and engineering applications. We address this issue by
exploiting parameter space reduction techniques in a multi-fidelity setting.

Gaussian processes (GP)~\cite{williams2006gaussian} have spread in
many fields as a reliable regression (GPR) method, especially for
optimization and inverse problems. Many extensions stemmed from the
original formulation, such as for kernel
methods~\cite{kanagawa2018gaussian}, and for big data and memory
limitations~\cite{lazaro2010sparse, liu2020gaussian}. On the other
hand the exploitation of multi-fidelity models had a huge impact in
the scientific computing community thanks to the possibility to
integrate simulations and data coming from different models and
sources. For an overview of different applications we suggest~\cite{kennedy2000predicting,
forrester2007multi, raissi2017inferring, bonfiglio2018multi,
bonfiglio2018improving, kramer2019multifidelity}.
A particularly promising nonlinear autoregressive multi-fidelity
Gaussian process regression (NARGP), was proposed in~\cite{perdikaris2017nonlinear}.
Recent advancements in the context of physics informed neural
networks~\cite{raissi2019physics} in a multi-fidelity setting for
function approximation and inverse PDE problems, can be found in~\cite{meng2020composite}.

These models achieve increased expressiveness with some kind of
nonlinear approach extending GP models to non-GP processes at the cost
of an additional computational load. In this direction, 
some works aim to obtain computationally efficient
heteroscedastic GP models using a variational inference
approach~\cite{lazaro2011variational}, or a nonlinear
transformation~\cite{snelson2004warped}. This approach is extended to
multi-fidelity models starting from the linear formulation presented by Kennedy and
O'Hagan~\cite{kennedy2000predicting} towards deep GP~\cite{damianou2013deep} and NARPG.

Classical low-fidelity models obtained by coarse grids or simplified
physical models still suffer the curse of dimensionality when used for
high-dimensional GP construction. Linear parameter space reduction
with Active Subspaces (AS)~\cite{constantine2015active}
can fight such curse using input-output couples obtained by
high-fidelity simulations.
Successful applications of parameter space reduction with active
subspaces can be found in many engineering fields: naval and nautical
problems~\cite{tezzele2018dimension}, shape
optimization~\cite{lukaczyk2014active, ghoreishi2019adaptive,
  demo2020asga, demo2021hull}, car aerodynamics
studies~\cite{othmer2016active}, inverse problems~\cite{constantine2016accelerating, nguyen2022dias}, 
cardiovascular studies coupled with intrusive model order
reduction~\cite{tezzele2018combined}, for the study of
high-dimensional parametric PDEs~\cite{o2022derivative}, and in CFD problems in a
data-driven setting~\cite{demo2019cras, tezzele2020enhancing}, among others.
New extensions of AS have also been developed in the recent years such
as AS for multivariate vector-valued
functions~\cite{zahm2020gradient}, a kernel approach for AS for scalar and vectorial
functions~\cite{romor2022kas}, a localization extension for both
regression and classification tasks~\cite{romor2021las}, and
sequential learning of active subspaces~\cite{wycoff2021sequential}. The multi-fidelity setting has been
used to find an active subspace given different fidelity models~\cite{lam2020multifidelity}.

Other nonlinear techniques for parameter space reduction include
manifold learning~\cite{holiday2019manifold, pozharskiy2020},
active manifolds~\cite{bridges2019active} and nonlinear level-set learning
(NLL)~\cite{zhang2019learning}. NLL adopts a Reversible Neural Networks
(RevNet) architecture to learn an effective parameter space
deformation to capture the geometry of the objective function
level-sets and parametrize them.

With this contribution we show how to integrate linear and nonlinear
parameter space dimensionality reduction within a multi-fidelity
regression scheme based of Gaussian processes to increase the accuracy
of high-dimensional response surfaces. The low-fidelity models are
built with AS or NLL and incorporated in the NARGP framework,
following the preliminary results obtained in~\cite{romor2020pamm}. An
extensive automotive test case is presented with different
configurations. 

This work is organized as follows: in Section~\ref{sec:multifidelity} we
introduce multi-fidelity Gaussian process regression starting from the
building block of a single fidelity up to the NARGP method; in
Section~\ref{sec:as} we focus on the parameter space reduction with
active subspaces and nonlinear level-set learning which are going
to be used to construct the low-fidelity models;
Section~\ref{sec:mfas} shows how to add the low-intrinsic
dimensionality bias into the NARGP framework, accompanied by
pseudocode; in Section~\ref{sec:results} we present the numerical
results of the proposed approach applied to two benchmark models, and to an
automotive application; finally
Section~\ref{sec:conclusions} draw the conclusions and some future
research lines.

\section{Multi-fidelity Gaussian process regression}
\label{sec:multifidelity}
In this section we are going to briefly recall the Gaussian process regression
(GPR) technique in order to better characterize the nonlinear autoregressive
multi-fidelity Gaussian process regression (NARGP) introduced
in~\cite{perdikaris2017nonlinear}. NARGP represents the main framework for our
proposed multi-fidelity method. We are going to consider the general
setting with multiple levels of fidelity.

\subsection{Gaussian process regression}

Gaussian process regression is a supervised technique to approximate unknown
functions given a finite set of input/output pairs $\mathcal{S} = \{ x_i, y_i
\}_{i=1}^N$. Let $f: \mathcal{X} \subset \mathbb{R}^m \to \mathbb{R}$ be the
scalar function of interest. The set $\mathcal{S}$ is generated through $f$ with
the following relation: $y_i = f(x_i)$, which are the noise-free observations.
We assigned a prior to $f$ with mean $m(\mathbf{x})$ and covariance function
$k(\mathbf{x}, \mathbf{x}^\prime ; \theta)$, that is $f (\mathbf{x}) \sim
\mathcal{GP} (m(\mathbf{x}), k(\mathbf{x}, \mathbf{x}^\prime ; \theta))$. The
prior expresses our beliefs about the function before looking at the observed
values. From now on we consider zero mean $\mathcal{GP}$, that is $m(\mathbf{x})
= \mathbf{0}$, and we define the covariance matrix as $\mathbf{K}_{i,  j} =
k(x_i, x_j ; \theta)$, with $\mathbf{K} \in \mathbb{R}^{N \times N}$. In order
to make predictions using the Gaussian process we still need to find the optimal
values for the hyper-parameters vector $\theta$ by maximizing the log
likelihood:
\begin{equation}
\log p(\mathbf{y} \vert \mathbf{x}, \theta) = - \frac{1}{2} \mathbf{y}^T \mathbf{K}^{-1} \mathbf{y} -\frac{1}{2}
\log \vert \mathbf{K} \vert - \frac{N}{2} \log 2 \pi.
\end{equation}
Let $\mathbf{x}_*$ be the test samples, and $\mathbf{K}_{N*} = k(\mathbf{x},
\mathbf{x}_* ; \theta)$ be the matrix of the covariances evaluated at all pairs
of training and test samples, and in a similar fashion $\mathbf{K}_{*N} =
k(\mathbf{x}_*, \mathbf{x} ; \theta)$, and $\mathbf{K}_{**} = k(\mathbf{x}_*,
\mathbf{x}_* ; \theta)$. By conditioning the joint Gaussian distribution on the
observed values we obtain the predictions $f_*$ by sampling the posterior as
\begin{equation}
f_* \vert \mathbf{x}_*, \mathbf{x}, \mathbf{y} \sim \mathcal{N} (\mathbf{K}_{*N} \mathbf{K}^{-1} \mathbf{y}, \mathbf{K}_{**} -
\mathbf{K}_{*N} \mathbf{K}^{-1} \mathbf{K}_{N*} ).
\end{equation}

\subsection{Nonlinear multi-fidelity Gaussian process regression}
\label{subsec:nargp}
In this section we briefly present the nonlinear autoregressive multi-fidelity
Gaussian process regression (NARGP) scheme~\cite{perdikaris2017nonlinear}. It extends the concepts present
in~\cite{kennedy2000predicting, le2014recursive} to nonlinear correlations
between the different available fidelities.

The procedure is purely data-driven. We start from the input/output pairs
corresponding to $p$ levels of increasing fidelity, that is
\begin{equation}
  \mathcal{S}_q = \{ x_i^q, y_i^q \}_{i=1}^{N_q} \subset
\mathcal{X} \times \mathbb{R} \subset \mathbb{R}^m \times \mathbb{R}, \qquad \text{ for } q \in \{1, \dots, p \},
\end{equation}
where $y_i^q = f_q (x_i^q)$. With $p$ we indicate the highest fidelity. We also
 assume that the design sets have a hierarchical structure:
\begin{equation}
  \pi(S_p)\subset\pi(S_{p-1})\subset\dots\subset\pi(S_1),
\end{equation}
where $\pi: \mathbb{R}^m\times\R\rightarrow \mathbb{R}^m$ is the projection onto
the first $m$ coordinates. Due to this hierarchy, when the fidelities of the
available datasets cannot be neatly assessed, it is resonable to consider the
cost needed to produce them as ordering criterion, see Remark~\ref{rmk:reversing
fidelities order}.

The NARGP formulation assigns a Gaussian process to each fidelity model $f_{q}$,
so they are completely defined by the mean field $m_{q}$, with the constant zero
field as prior, and by their kernel $k_{q}$, as follows:
\begin{equation}
  y_q(\bar{x})-\epsilon \sim \mathcal{GP}(f_q(\bar{x})) \vert  m_q(\bar{x}), k_q(\theta_q))\quad \forall q \in \{1, \dots, p \}\; ,
\end{equation}
where $\epsilon\sim\mathcal{N}(0, \sigma^{2})$ is a noise term and
\begin{equation}
  \bar{x}:=\begin{cases}
    (\mathbf{x},f_{q-1}(\mathbf{x}))\in \mathbb{R}^m\times \mathbb{R},&q>1\\
    \mathbf{x}\in \mathbb{R}^m,&q=1
  \end{cases}.
\end{equation}

The definition of the kernel $k_q(\theta_q)$ implements the auto-regressive
characteristic of the method since it depends on the previous fidelity model
$f_{q-1}$:
\begin{equation}
  \label{eq:autoregressive_kernel}
  k_q((x, f_{q-1}(x)), (x', f_{q-1}(x'));\theta_q) = k^{\rho}_q(x,
  x';\theta_q^{\rho})\cdot k^f_q(f_{q-1}(x), f_{q-1}(x');\theta^f_q) +
  k^{\delta}_{q}(x, x';\theta^{\delta}_q)\; .
\end{equation}
The hyper-parameters to be tuned are represented by
$\theta_q\equiv(\theta_{q}^{\rho}, \theta_{q}^{f}, \theta_{q}^{\delta})$ and are
associated respectively to the multiplicative kernel $k_{q}^{\rho}$, the
auto-regressive kernel $k^{f}$, and the kernel $k_{q}^{\delta}$, which
corresponds to the non auto-regressive part in the sum of
Equation~\ref{eq:autoregressive_kernel}. For our applications we employ the
radial basis function kernel with automatic relevance determination
(RBF-ARD)~\cite{williams2006gaussian}, but there are other possible choices.

The presence of the multiplicative kernel $k_{q}^{\rho}$ allows nonlinear
interdependencies between subsequent fidelities to be modelled, surpassing a
linear auto-regressive multi-fidelity scheme. The latent manifold that relates
the inputs, the lower fidelity posterior and the high-fidelity posterior is in
this case nonlinear~\cite{perdikaris2017nonlinear}.

We use the notation$(\mathbf{x}, y_{q-1}(\mathbf{x}))$ for
the training set and $\mathbf{x}_{*}$ for the new input. So in order to evaluate
the predictive mean and variance for a new input $\mathbf{x}_{*}$ we have to
integrate the posterior $p(f_{q}(\mathbf{x}_{*}) \vert f_{q-1}, \mathbf{x}_{*},
\mathbf{x}_{q}, y_{q})$ defined as
\begin{align}
&f_{q}(\mathbf{x}_{*} \vert f_{q-1}, \mathbf{x}_{*}, \mathbf{x}_{q}, y_{q})\sim\mathcal{N} (\mathbf{K}^{q}_{*N} (\mathbf{K}^{q})^{-1} y_{q}, \mathbf{K}^{q}_{**} -
\mathbf{K}^{q}_{*N} (\mathbf{K}^{q})^{-1} \mathbf{K}^{q}_{N*} ),\\
&\mathbf{K}_{*N}^{q}=k_{q}((\mathbf{x}_{*},f_{q-1}(\mathbf{x}_{*}), (\mathbf{x}_{q-1}, y_{q-1});\theta_q),\\
&\mathbf{K}_{N*}^{q}=k_{q}((\mathbf{x}_{q-1}, y_{q-1}), (\mathbf{x}_{*},f_{q-1}(\mathbf{x}_{*});\theta_q),\\
&\mathbf{K}^{q}=k_{q}((\mathbf{x}_{q-1}, y_{q-1}), (\mathbf{x}_{q-1}, y_{q-1});\theta_q),
\end{align}
over the Gaussian distribution of the prediction at the previous level
$f_{q-1}(\mathbf{x}_{*})\sim\mathcal{N}(m_{q-1}(\mathbf{x}_{*}),k_{q-1}(\mathbf{x}_{*}))$.
Apart from the first level of fidelity $q=1$ the posterior probability
distribution given the previous fidelity models is no longer Gaussian.
So, in practice, the following integral is approximated with recursive Monte Carlo at
each fidelity level, for all $q\in\{2,\dots, p\}$,
\begin{align}
p(f_{q}^{\text{post}}(\mathbf{x}_{*}))&:=p(f_{q}(\mathbf{x}_{*}) \vert
                                        f_{q-1}, \mathbf{x}_{*},
                                        \mathbf{x}_q,
                                        y_q) = \nonumber \\ &=\int_{\mathcal{X}}
                                        p(f_{q}(\mathbf{x}_{*}) \vert s, \mathbf{x}_{*}, \mathbf{x}_{q}, y_{q}) d\mathcal{L}_{f_{q-1}^{\text{post}}(\mathbf{x}_{*})}(s)\\
p(f_{1}^{\text{post}}(\mathbf{x}_{*}))&:=p(f_{1}(\mathbf{x}_{*}) \vert
                                        \mathbf{x}_{*}, \mathbf{x}_1, y_1)\sim\mathcal{N}(m_1(\mathbf{x}_{*}),k_1(\mathbf{x}_{*})),
\end{align}
where $\mathcal{L}_{f_{q-1}^{\text{post}}(\mathbf{x}_{*})}$ is the probability
law of
$f_{q-1}^{\text{post}}(\mathbf{x}_{*})\sim\mathcal{N}(m_{q-1}(\mathbf{x}_{*}),k_{q-1}(\mathbf{x}_{*}))$.
In the applications we always use $200$ to $10000$ Monte Carlo samples, since the results
do not vary much increasing them for our test cases.

The hyper-parameters $\theta_q$ are optimized (non recursively) with maximum
log-likelihood estimation for each GP model $\mathcal{GP}(f_q \vert \mathbf{0},
k^q(\theta_q))$, for all $q\in\{1,\dots,p\}$,
\begin{equation*}
  \argmin_{\theta_q}-\log{p(f_q(\mathbf{x}_q) \vert \mathbf{x}_q, y_q,
    y_{q-1},\theta_q)}\propto \frac{1}{2}\log{\vert K^q(\theta_q)\vert
  }+\frac{1}{2}y_q^{T}(K^q(\theta_q))^{-1}y_q\; ,
\end{equation*}
this is why a hierarchical dataset is needed. The hyperparameters tuning is
achieved maximizing the log-likelihood with the gradient descent optimizer
L-BFGD in GPy~\cite{gpy}. For some test cases, the training procedure is subject
to relevant perturbations relative to the number of restarts, this is especially
true in higher dimensions of the parameter space.

\section{Parameter space reduction}
\label{sec:as}
Our aim is testing multi-fidelity Gaussian process regression models to
approximate objective functions which depend on inputs/parameters sampled from a
high-dimensional space. Low-fidelity models relying on a physics-based or
numerical model reduction --- for example a coarse discretization or a more
specific numerical model order reduction --- still suffer from the high
dimensionality of the input space. In our approach we try to tackle these problematics by searching for a surrogate
(low-fidelity) model accounting for the complex correlations among the input
parameters that concur to the output of interest. With this purpose in mind, in
this section we are going to briefly present the active subspaces (AS)~\cite{constantine2015active}, and the nonlinear level-set learning
(NLL) method~\cite{zhang2019learning} for parameter space reduction in order to
design response surfaces with Gaussian process regression.

\subsection{Active subspaces}
\label{subsec:as}
Let $\mathbf{X}$ be an absolutely continuous random variable with
probability density $\rho$, such that
$\text{supp}(\rho)=\mathcal{X}\subset\mathbb{R}^{m}$. The variable
$\mathbf{X}$ represents the inputs, and $m$
denotes the dimension of the input parameter space. With simple Monte
Carlo we can approximate the uncentered covariance matrix of the gradients of
the function of interest as
\begin{equation}
\mathbb{E}_{\rho}[\nabla_{\mathbf{x}} f (\nabla_{\mathbf{x}}
f)^{T}]\approx\frac{1}{N}\sum_{i=1}^{N} \nabla_{\mathbf{x}} f(\mathbf{X}_{i})
(\nabla_{\mathbf{x}} f(\mathbf{X}_{i}))^{T} ,
\end{equation}
where $N$ denotes the number of samples. We are looking for the highest spectral gap
$\lambda_{r}-\lambda_{r+1}$ in the sequence of ordered eigenvalues of the
approximated correlation matrix.
The active subspace is the eigenspace corresponding to the first $r$ eigenvalues
$\lambda_{1},\dots,\lambda_{r}$ and it is denoted with the matrix
\RA{$\hat{W_r}\in\mathcal{M}(m\times r)$} whose columns are the corresponding $r$ active
eigenvectors. The inactive subspace is defined as the span of the
remaining eigenvectors. On it $f$ is almost flat on average, so we can
safely discard such component without compromising too much the
accuracy. We can thus build a response surface $\mathcal{R}$ using a Gaussian
process regression trained with $N_{\text{train}}$ pairs \RA{$\{\hat{W}_r^{T}\mathbf{x}_{i},
y_{i}\}_{i=1}^{N_{\text{train}}}$} of active inputs and outputs.

The mean square regression error is bounded a priori~\cite{constantine2015active} by
\RA{\begin{equation}
    \mathbb{E}_{\rho}\left[(f(\mathbf{X})-\mathcal{R}(\hat{W}_r^{T}\mathbf{X}))^2\right]
    \leq C_1(1+N^{-1/2})^2\left(
      \epsilon(\lambda_1+\dots+\lambda_r)^{1/2}+(\lambda_{r+1}+\dots +
      \lambda_m)^{1/2}\right)^2+C_2\delta,
\end{equation}}
where $C_{1}$ and $C_{2}$ are constants, $\epsilon$ quantifies the error in the
approximation of the true active subspace \RA{$W_r$ with $\hat{W}_r$} obtained
from the Monte Carlo approximation, and $C_{2}\delta$ is a bound on the mean
squared error of the Gaussian process regression over the active subspace:
\RA{\begin{equation}
  \mathbb{E}_{\rho_{\hat{W}_r\X|\hat{W}_{m-r}\X}}\left[\left(\overline{\mathbb{E}_{\rho}\left[
          f(\X) | \sigma(\hat{W}_r^T\X)
        \right]}-\mathcal{R}(\hat{W}_r^{T} \mathbf{X})\right)^2\right]\leq C_2\delta,
\end{equation}}
where \RA{$\rho_{\hat{W}_r\X|\hat{W}_{m-r}\X}$} is the probability of the active
variables conditioned on the inactive ones, and
\RA{$\overline{\mathbb{E}_{\rho}\left[ f(\X) | \sigma(\hat{W}_r^T\X) \right]}$} is
the random variable $f(\X)$ conditioned on the $\sigma$-algebra generated by
\RA{$\hat{W}_r^T \X$} and approximated with the Monte Carlo method.

\subsection{Nonlinear level-set learning method}
\label{subsec:nll}
This method seeks a bijective nonlinear transformation
$g_{\text{NLL}}:\mathcal{X}\rightarrow \Tilde{\mathcal{X}}\subset\mathbb{R}^{m}$
to capture the geometry of level sets and parametrize them in a low-dimensional
space. To this end in~\cite{zhang2019learning} they employ reversible networks
(RevNets)~\cite{chang2018reversible} to learn the transformation
$g_{\text{NLL}}$. The designed loss function uses samples of the gradients of
the target function to encourage the transformed function to be sensitive to
only a few active coordinates.


To construct the RevNet, the following
architecture~\cite{haber2017stable}, which is reversible by definition, is employed:
\begin{equation}
\left \{
\begin{array}{l}
  u_{n+1} = u_n + h K^T_{n, 1} \sigma (K_{n, 1} v_n + b_{n, 1}) \,,\\
  v_{n+1} = v_n - h K^T_{n, 2} \sigma (K_{n, 2} u_n + b_{n, 2}) \,,
\end{array} \right .
	\qquad \mathrm{for}\;\; n = 0, 1, \dots , N-1,
\end{equation}
where $u$ and $v$ are partitions of the states, $h$ is a scalar time step, the
matrices $K$ contain the weights, $b$ represent the biases, and $\sigma$ is the
activation function. We remark that the original coordinates and the transformed
ones are split in two in $u$ and $v$.

\section{Multi-fidelity data fusion with active subspaces}
\label{sec:mfas}
Our study regards the design of a nonlinear autoregressive multi-fidelity
Gaussian process regression (NARGP)~\cite{perdikaris2017nonlinear} with two
fidelities: the high-fidelity corresponds to a relatively accurate and costly
model, for example a numerical model which requires computationally intensive
simulations to obtain a scalar output for each parameter sample; and the
low-fidelity level which comes from a
response surface built through a parameter
space reduction technique --- here we focus on active subspaces but little
modifications are required in order to use NLL as we are going to
show. We consider models with high-dimensional input space but with a
low intrinsic dimensionality. This setting characterizes many
industrial applications~\cite{othmer2016active,
  tezzele2018dimension, lukaczyk2014active}.

In fact, the inductive biases we impose come mainly from two sources: the kernel
of the Gaussian process (lengthscale, noise, regularity of the stochastic process) and the low-fidelity
intrinsic dimensionality assumption (presence of a dominant linear or nonlinear
active subspace). The key feature of the method is the imposition of the latter
on the multi-fidelity model design: we expect that a hint towards the
presence of an active subspace will be transferred from the low-fidelity to the
high-fidelity level through the discovery of nonlinear correlations between the
low-fidelity predictions, and the high-fidelity inputs/outputs dataset.
In this way, the accuracy should increase in the data-scarcity regime, i.e. when the
number of high-fidelity samples are not enough to obtain an accurate
single-fidelity regression.
The overhead with respect to
the original procedure~\cite{perdikaris2017nonlinear} is the evaluation of the
active subspace from the high-fidelity inputs and the training of the whole
multi-fidelity model; this costs are usually negligible as shown in section~\ref{sec:automotive}.

\begin{figure}[ht!]
  \centering
  \includegraphics[width=1.\textwidth]{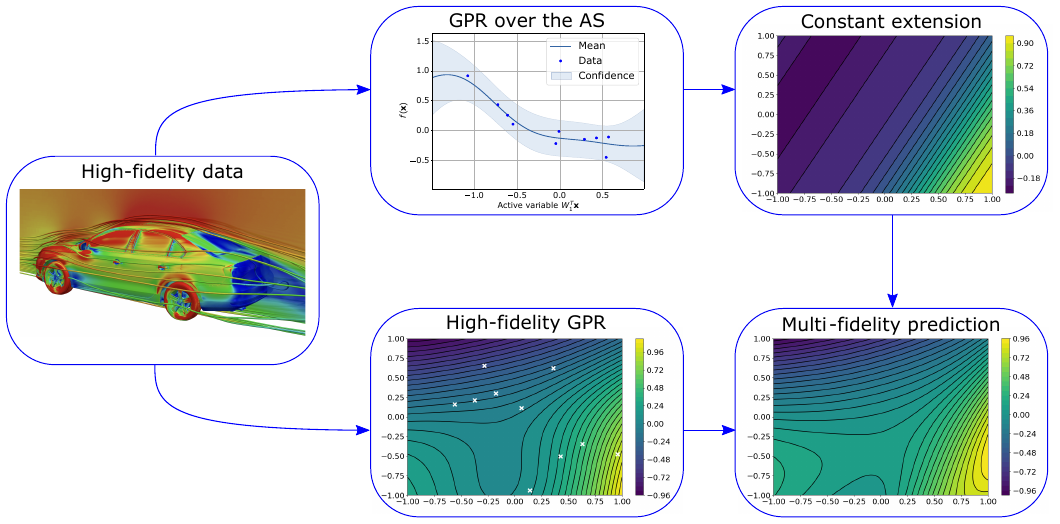}
  \caption{Illustrative scheme of the NARGP-AS method. Starting from $10$
    high-fidelity data (depicted with blue dots and white crosses) we construct as
    low-fidelity model a response surface which is constant along the inactive
    subspace.}
  \label{fig:nargpas}
\end{figure}

In Figure~\ref{fig:nargpas} we present an illustrative scheme
of the proposed NARGP-AS method; the underling objective function is an
hyperbolic paraboloid $f:[0,1]^2\subset\mathbb{R}^2\rightarrow\mathbb{R},\ f(x_1,
x_2)=x_1^2-x^2_2$ and is shown only for the purpose of
representing the procedure more clearly. The high-fidelity flow field belongs to
the automotive application of section~\ref{sec:automotive}.

For clarity we will use the letters $H$ and $L$ as labels for the high-fidelity
and low-fidelity models respectively, instead of the fidelity levels $q=1$ and $q=2$.
Changing, as just described, the notations of subsection~\ref{subsec:nargp}, we consider
\begin{align*}
  &S^{L} = \{x_i^{L}, y_i^{L}\}^{N_L}_{i=1}\subset \mathbb{R}^{m}\times\mathbb{R},\\
  &S^{H} = \{x_i^{H}, y_i^{H}\}^{N_H}_{i=1}\subset \mathbb{R}^{m}\times\mathbb{R},
\end{align*}
and additionally $\{ dy^H_i\}_{i=1}^{N_2}\subset \mathbb{R}^{m}$ the gradients
corresponding to the high-fidelity dataset $S^{H}= \{x_i^{H}, y_i^{H}\}^{N_H}_{i=1}$: in principle, the gradients
can be directly obtained from the model of interest (with adjoint methods in case of PDE
models for example) or approximated from the input output pairs $S^{H}$. \RB{For the influence of the gradients' approximation on the regression error, see~\cite{constantine2015active}. For our test cases we employ the exact gradients when available (for the benchmarks in sections~\ref{subsubsec:piston} and~\ref{subsubsec:ebola}) or we approximate them from the high-fidelity GPR (in the test cases in sections~\ref{subsubsec:Jetta-6}, \ref{subsubsec:Jetta-12-RANS} and~\ref{subsubsec:Jetta-12-RANS-DDES}). This does not results in additional costs, since the HF GPR is needed in the NARGP-AS procedure.}

The high-fidelity dataset $S^H$ (and the corresponding gradients) represents by
itself all the necessary ingredients: $S^L$ is built from $S^H$ through a
response surface on the linear or nonlinear active subspace. For this purpose
the dataset $S^H$ is employed, with the corresponding gradients, to find an
active subspace \RA{$\hat{W}_r$} or train a RevNet, as described in subsections
\ref{subsec:as} and \ref{subsec:nll}.

Then, since $S^H\subset
S^L$, we write
\begin{equation}
  S^L\setminus S^H =\{x_i^L, y_i^L\}_{i=1}^{N_L}\setminus \{x_i^H, y_i^H\}_{i=1}^{N_H} = \{\Tilde{x}_i^L, \Tilde{y}_i^L\}^{N_L-N_H}_{i=1}.
\end{equation}
The additional low-fidelity inputs $\{\Tilde{x}_i^1\}^{N_L-N_H}_{i=1}$ are sampled
independently from the inputs' probability distribution, while the additional low-fidelity
outputs $\{\Tilde{y}_i^L\}^{N_L-N_H}_{i=1}$ are the predictions associated to the
active components of the additional low-fidelity inputs
\RA{$\{\hat{W}_r^{T}\Tilde{x}_i^L\}^{N_L-N_H}_{i=1}$}, obtained from the response
suface trained on
\RA{$$\{\hat{W}_r^{T}x_i^{H},
y_i^{H}\}^{N_H}_{i=1}\subset\mathbb{R}^{r}\times\mathbb{R}.$$}
The response surface is
trained as a Gaussian process regression as described in
subsection~\ref{subsec:as}. The procedure is synthetically reviewed through
Algorithm~\ref{algo:mfas_dependent}. \RA{The number of low-fidelity samples is chosen until a good approximation of the low-fidelity response surface is obtained. As it is experimentally shown in Figure~\ref{fig:r2_scores_jetta6}, additional low-fidelity samples do not improve the accuracy of the multi-fidelity model afterwards.}

\begin{remark}[Nonlinear level-set learning as LF model]
    If NLL is employed to build the low-fidelity level, only the first step of
    Algorithm~\ref{algo:mfas_dependent} is changed. For our applications, the
    GPR designed with NLL has dimension one.
\end{remark}

\begin{algorithm}[ht]
    \caption{NARGP-AS response surface design algorithm.}
    \label{algo:mfas_dependent}

    \SetKwInOut{Input}{input}\SetKwInOut{Output}{output}

    \Input{training high-fidelity inputs, outputs, gradients triplets $\{(\x^{H}_{i},
        y^{H}_{i}, dy^{H}_{i})\}_{i=1}^{N_H}\subset \R^{m}\times \R\times
        \R^{m}$,\\
    low-fidelity inputs $\{\x^{L}_{i}\}_{i=1}^{N_L}\subset \R^{m}$,\\
    } \BlankLine \Output{multi-fidelity model, $\left((f_{H}|x^{H}_{i}, y^{H}_{i}),\ (f_{L}|x^{L}_{i})\right)\sim
                \left(\mathcal{GP}(f_{H}|m_{H}, k_{H}),
                \mathcal{GP}(f_{L}|m_{L}, k_{L})\right)$}

    \BlankLine Compute the active subspace $\hat{W}_{r}$ with the high-fidelity
    gradients $\{dy_{i}^{H}\}_{i=1}^{N_H}$,\\
    Build the response surface $\mathcal{R}(\hat{W}_{r}\X)$ with
    a GP regression from $\{(\hat{W}_{r}\x^{H}_{i}, y^{H}_{i})\}_{i=1}^{N_H}$,\\
    Predict the low-fidelity outputs $\{y^{L}_{i}\}_{i=1}^{N_L}$ at
    $\{\x^{L}_{i}\}_{i=1}^{N_L}$ and the training high-fidelity inputs $\{y^{H}_{i}\}_{i=1}^{N_H}$ at $\{\x^{H}_{i}\}_{i=1}^{N_H}$ with the
    response surface,\\
    Train the multi-fidelity model at the low-fidelity level $f_{L}$ with the
    training dataset $\{(x^{L}_{i}, y^{L}_{i})\}_{i=1}^{N_L}\cup\{(x^{H}_{i},
    y^{H}_{i})\}_{i=1}^{N_H}$,\\
    Train the multi-fidelity model at the high-fidelity level $f_{H}$ with the
    training dataset $\{((x^{H}_{i}, y^{H}_{i}),
    y^{H}_{i})\}_{i=1}^{N_H}$

    \BlankLine
  \end{algorithm}

\begin{remark}[Markov property]
    Theoretically the observations $\{y^{q}_{i}\}$ should be noiseless for each
level of fidelity $q$ in order to preserve the Markov
property~\cite{perdikaris2017nonlinear}. However, in practice, it could be beneficial
in some applications to add noise at each fidelity level, or constraint the noise levels from below in order to avoid overfitting.
\label{rmk:Markov}
\end{remark}

\section{Numerical results}
\label{sec:results}
In this section we are going to present the results obtained with the NARGP-AS
and the NARGP-NLL method over two benchmark test problems (Piston~\ref{subsubsec:piston} and Ebola~\ref{subsubsec:ebola} models), and over a more
complex car aerodynamics problem (Jetta-6~\ref{subsubsec:Jetta-6}, Jetta-12-RANS~\ref{subsubsec:Jetta-12-RANS}, Jetta-12-DDES~\ref{subsubsec:Jetta-12-RANS-DDES}). The library employed to implement the NARGP
model is Emukit~\cite{emukit2019} while for the active subspace and NLL response
surface design we used the open source Python package\footnote{Available at
\url{https://github.com/mathLab/ATHENA}.} called ATHENA~\cite{romor2020athena},
and GPy~\cite{gpy}.

\RB{The computational times of the prediction and training of the
  NARGP-AS method are reported in Table~\ref{tab: computational
    times}. In particular, it is shown how the number of HF test samples and of Monte Carlo (MC) samples affect the MF prediction times. The training costs are mainly affected by the number of restarts of the optimization with L-BFGD, instead.}
\begin{table}[htp!]
  \centering
  \caption{\RB{Computational times of the training of the multi-fidelity models and evaluation of the predictions.}}
  \footnotesize
  \begin{tabular}{ l  c  c  c  c  c  c  c }
      \hline
      \hline
      \multirow{2}{*}{Test Case} & \# HF Training & \multirow{2}{*}{Training} & \multirow{2}{*}{Restarts} & MC & \# HF Test & MF & HF\\
    & Samples & & & Samples & Samples & Prediction & Prediction \\
    \hline
    \hline
    \rowcolor{Gray}
      Piston model~\ref{subsubsec:piston}  & 150  & 24 [s] & 10 & 100 & 10000  & 10 [s] & 0.123 [s]   \\
      \hline
      Ebola model~\ref{subsubsec:ebola}   & 150  & 21 [s] & 10 & 100 & 10000  & 10 [s] & 0.450 [s]     \\
    \hline
    \rowcolor{Gray}
      Jetta-6~\ref{subsubsec:Jetta-6} & 76 & 229 [s] & 150  & 100 & 25  & 0.056 [s] & 0.0006 [s]      \\
      \hline
      Jetta-12-RANS~\ref{subsubsec:Jetta-12-RANS}  & 185 & 50 [s] & 10 & 10000 & 51  & 14.2 [s] & 0.0006 [s]      \\
    \hline
    \rowcolor{Gray}
      Jetta-12-DDES~\ref{subsubsec:Jetta-12-RANS-DDES}  & 65 & 20 [s] & 10 & 1000 & 50  & 0.02 [s] & 0.0003 [s]     \\
      \hline
      \hline 
  \end{tabular}
  \label{tab: computational times}
\end{table}

\subsection{Benchmark test problems}
The first benchmark test problem presents a $7$-dimensional input
parameter space and the quantity of interest is the time a cylindrical
piston takes to complete a cycle\footnote{The piston dataset was
taken from \url{https://github.com/paulcon/active_subspaces}.}. The
second one is a $8$-dimensional model for the spread of Ebola in
Western Africa~\cite{diaz2018modified}. These tests have been chosen because of the presence of an active subspace
and they indeed present a low intrinsic dimensionality. The sufficient
summary plot is plotted for both the cases together with a
one-dimensional Gaussian process regression built over the AS. We also show
the correlation between the low-fidelity level and the high-fidelity level of the
multi-fidelity model. We compare the performance of the different
fidelities looking at the corresponding $R^2$ scores. \RA{This score
  is chosen to show how the obtained regressions compare with respect
  to a constant predictor equal to the function average ($R^2 = 0$)}. With LF we
denote the low-fidelity model represented by a GP regression on the
low-fidelity input/output couples, with HF the high-fidelity model
represented by a GP regression built on the full space, and with MF
the proposed multi-fidelity model. The number of low-fidelity samples
is kept fixed at $200$ for both test cases, while we study the accuracy varying the number of high-fidelity training samples used. For both the benchmark problems the models were tested over a dataset comprising $10000$ samples, selected with Latin hypercube sampling (LHS). The
nonlinear autoregressive fidelity fusion approach achieves better
performance with a consistent \RA{increase} in the $R^2$ score.

\subsubsection{The piston model}
\label{subsubsec:piston}

For this model the scalar target function of interest represents the
time it takes the piston to complete a cycle, depending on a
$7$-dimensional parameters vector. For its evaluation
a nonlinear function has to be computed. The input
parameters are uniformly distributed. For a detailed description of
the parameters' ranges the reader can refer
to~\cite{constantine2017global}.
The algebraic cylindrical piston model appeared as a test for statistical
screening in~\cite{ben2007modeling}, while
in~\cite{constantine2017global} they describe an active subspaces analysis.

\begin{figure}[ht!]
  \centering
  \includegraphics[width=.49\textwidth]{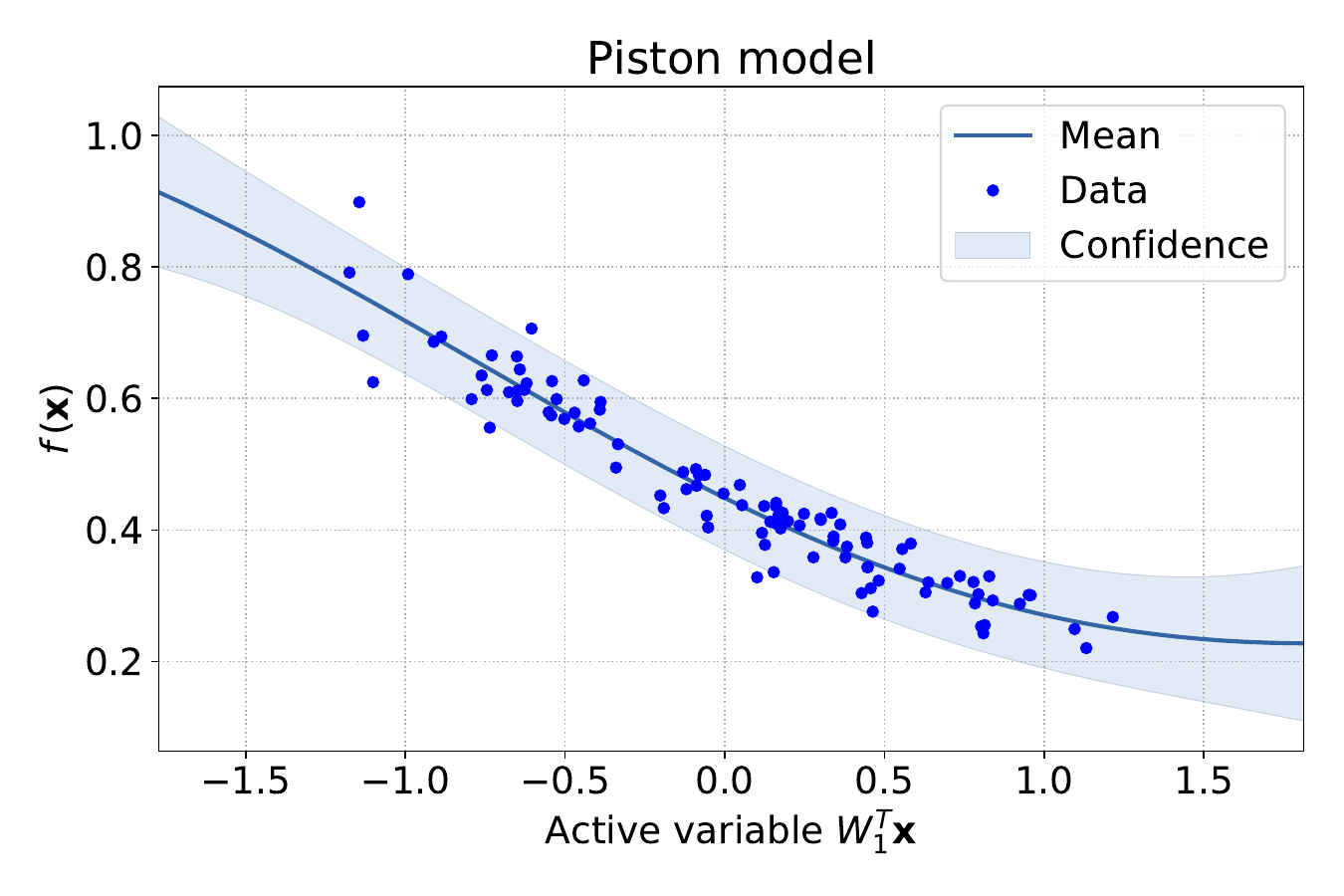}\hfill
  \includegraphics[width=.49\textwidth]{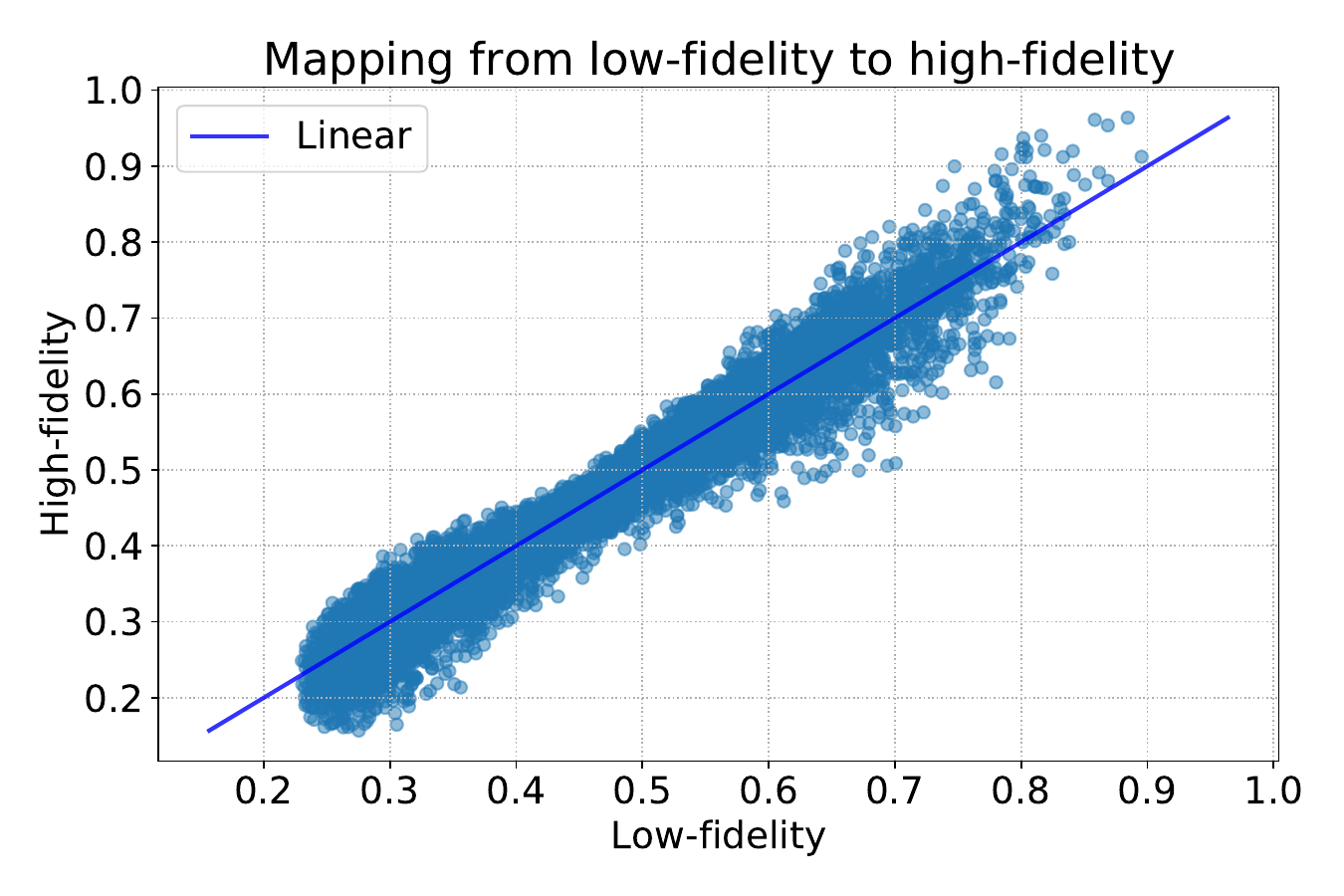}
  \caption{Left: sufficient summary plot of the surrogate model built with
    active subspaces. 100 samples were used to build the AS surrogate model
    shown. Right: correlations among the low-fidelity level and the high-fidelity
    level of the multi-fidelity model, evaluated at the $10000$ test samples.}
  \label{fig:piston_correlations}
\end{figure}

\begin{figure}[ht!]
  \centering
  \includegraphics[width=1.\textwidth]{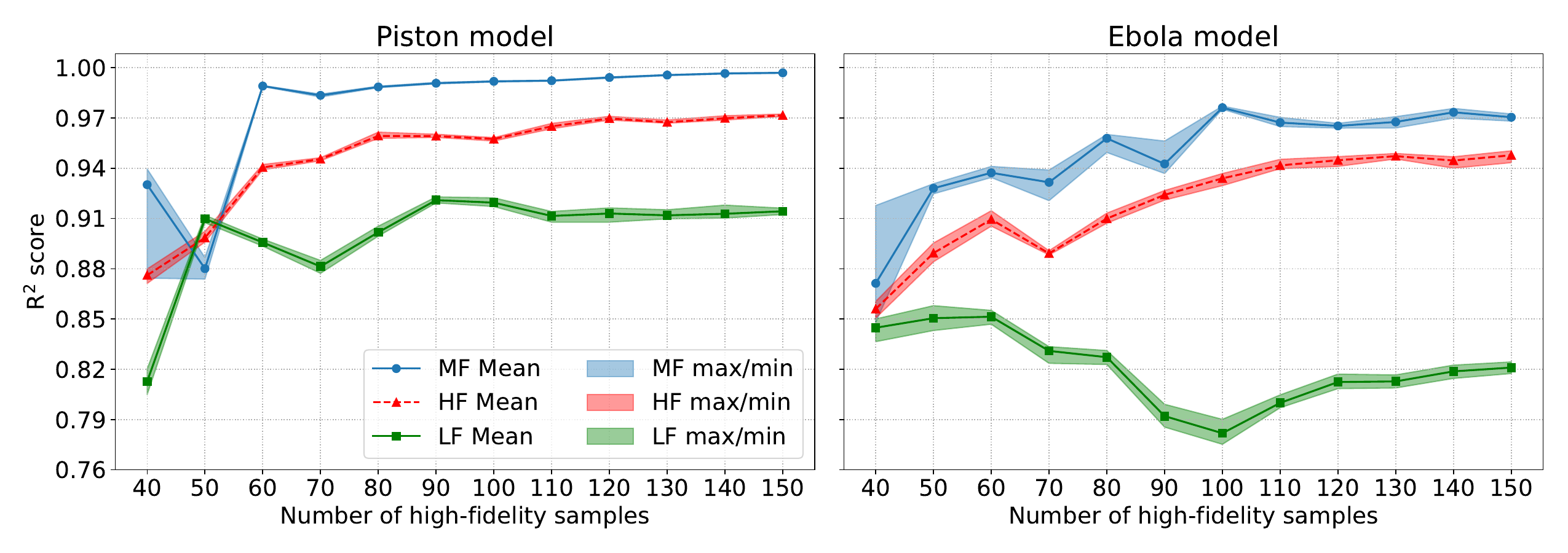}
  \caption{$R^{2}$ score of the posterior of the multi-fidelity (MF),
    high-fidelity (HF) and low-fidelity (LF) models against the number of
    high-fidelity samples used to find the active subspace and build the
    Gaussian process regressions of the MF, HF, LF models. The 10000 test
    samples are distributed with Latin hypercube sampling (LHS). In the left
    panel the results for the piston model, while on the right the Ebola spread
    model.}
  \label{fig:piston_ebola_mf}
\end{figure}

From the sufficient summary plot reported in the left panel
of Figure~\ref{fig:piston_correlations} we can conclude that a
one-dimensional active subspace is able to describe the input-output
dependency with a sufficient accuracy. This is also supported by the
GPR built over the AS. Moreover, the ordered eigenvalues of the
covariance matrix of the gradients exhibit a spectral gap between the
first and the second eigenvalue. In the right panel of
Figure~\ref{fig:piston_correlations} we present the correlation
between the high- and low-fidelity of the NARGP model.

Figure~\ref{fig:piston_ebola_mf} shows on the left the mean $R^2$ scores of the
MF model built as described in Section~\ref{sec:mfas} varying the
number of high-fidelity data. This is done over $10$ training
restarts of the MF, LF and HF models: moreover each GPR training is restarted
$10$ times for the HF and LF models and $20$ times for the MF model at each
fidelity level, inside the GPy package. We show also the minimum and maximum
$R^2$ scores over the outer $10$ training restarts to show the stability of the
procedure. When we have a scarce amount of data the models are not so
robust as we can see in the left part of the plot for $50$ and $60$
high-fidelity samples. After that point we have very stable results
which account for a relative gain in the $3$--$5\%$ range with respect to the
high-fidelity regression.

\subsubsection{Modified SEIR model for Ebola}
\label{subsubsec:ebola}

Now we consider the modified SEIR model for the spread of Ebola in
Liberia, presented in~\cite{diaz2018modified}, which depends on $8$ parameters. As
scalar output of interest we take the basic reproduction number
$R_0$. It can be computed with the following formula:
\begin{equation}
  \label{eq:Ebola}
  R_0 =\frac{\beta_1 +\frac{\beta_2\rho_1 \gamma_1}{\omega} +
    \frac{\beta_3}{\gamma_2} \psi}{\gamma_1+ \psi},
\end{equation}
with parameters range taken from~\cite{diaz2018modified}, where they
conducted a global sensitivity analysis with AS. For a kernel-based
active subspaces comparison the reader can refer to~\cite{romor2022kas}.

\begin{figure}[ht!]
  \centering
  \includegraphics[width=.49\textwidth]{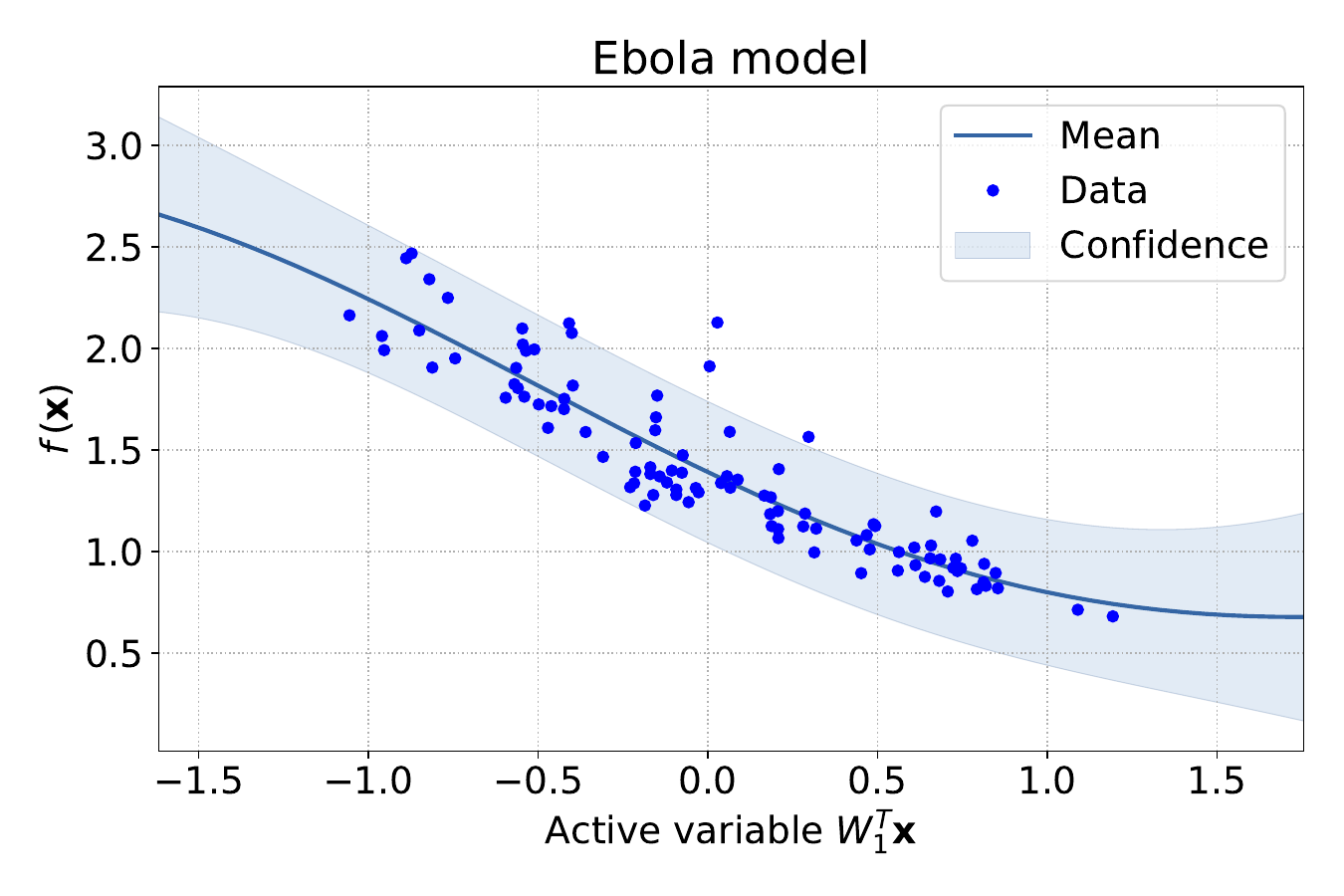}\hfill
  \includegraphics[width=.49\textwidth]{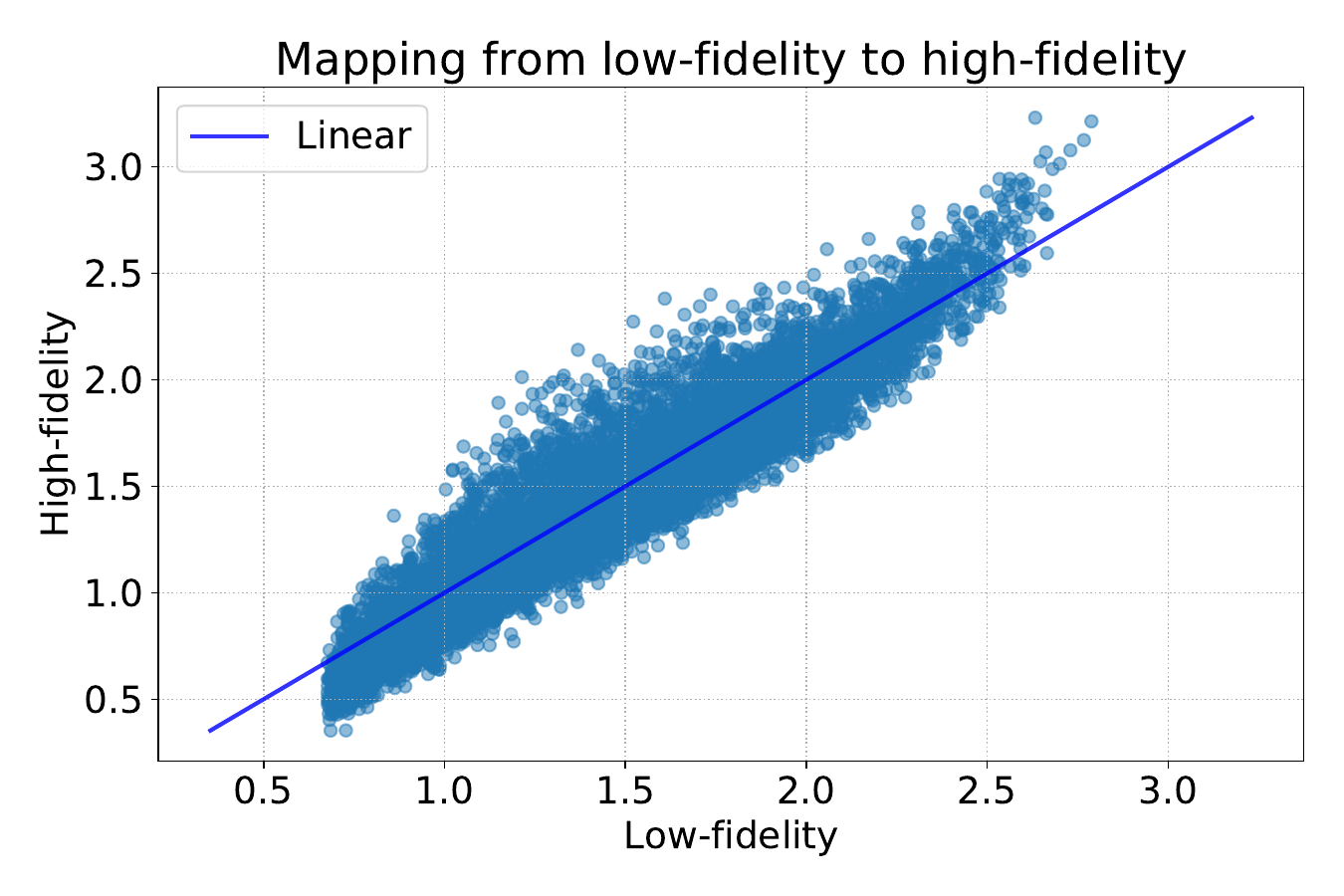}
  \caption{Left: sufficient summary plot of the Ebola model, 100 samples were
    used to build the AS surrogate model shown. Right: correlations among the
    low-fidelity level and the high-fidelity level of the multi-fidelity model, evaluated at the $10000$ test samples.}
  \label{fig:Ebola_correlations}
\end{figure}

In this case a one-dimensional Gaussian process response surface is
not able to achieve the same good accuracy of the previous case, as
can be seen in the left panel of
Figure~\ref{fig:Ebola_correlations}. This is also confirmed by the
correlation between the low- and high-fidelity levels of the NARGP, depicted in the
right panel of Figure~\ref{fig:Ebola_correlations}. The corresponding
$R^2$ scores in the right panel of Figure~\ref{fig:piston_ebola_mf}
reflect this behaviour of worse performance with respect to the piston
test case, where better correlations among the fidelities were
identified. The relative gain is in the $3$--$4\%$ range with respect
to the high-fidelity regression.

\subsection{Automotive application}
\label{sec:automotive}
\begin{figure}[ht!]
  \centering
  \includegraphics[width=.8\textwidth, trim=0 0 0 0, clip]{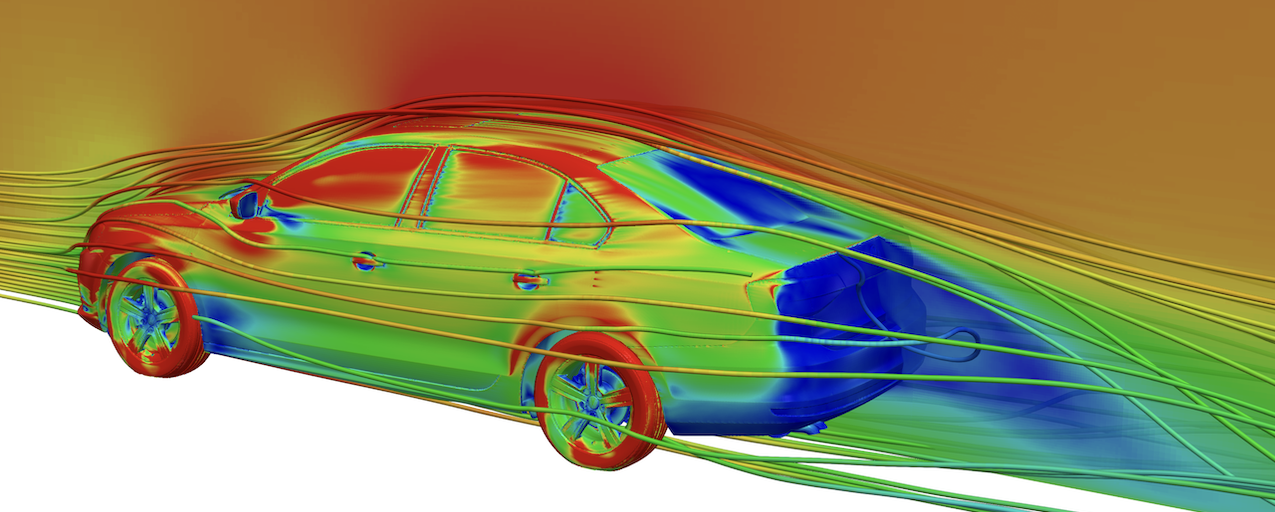}
  \caption{Visualization of the averaged flow field around the Jetta.}
  \label{fig:jetta_12_flow}
\end{figure}

Two different test cases from the world of automotive aerodynamics are
investigated in order to demonstrate the applicability of the presented method
to real-life problems. The first one (named hereafter {\bf Jetta-6}) is taken
from~\cite{mrosek2019}, where it is described in detail. It consists of a
$6$-dimensional geometric parameterization of the Volkswagen Jetta VI. The
parameters (see Table~\ref{tab:jetta_6_params}) were generated by free-form
deformation and focus on the rear part of the car. A Latin Hypercube with $101$
samples was created, and the aerodynamic flow fields were computed with
OpenFOAM~\cite{of} via Delayed Detached Eddy Simulations (DDES). An illustrative
example can be seen in Figure~\ref{fig:jetta_12_flow}. The physical simulation
time was four seconds, and the fields were averaged over the last two seconds
before integrating them over the vehicle surface to obtain the drag coefficient
$c_D$. With mesh sizes being of the order of $100$M cells, each variant required
about 23,000 CPU-core-hours.

\begin{table}[htb]
\centering
\caption{Parameters' description of the Jetta-6 test
case~\cite{mrosek2019}.\label{tab:jetta_6_params}}
\begin{tabular}{ c l r r }
\hline
\hline
Parameter & Description & Lower bound & Upper bound \\
\hline
\hline
\rowcolor{Gray}
$\mupar_1$  & Rear roof lowering  & 0 mm & 50 mm \\
$\mupar_2$  & Trunk height  & -30 mm & 30 mm \\
\rowcolor{Gray}
$\mupar_3$  & Trunk length  & -50 mm & 100 mm \\
$\mupar_4$  & Rear lateral tapering  & -60 mm & 50 mm \\
\rowcolor{Gray}
$\mupar_5$  & Rear end edge position  & -70 mm & 30 mm \\
$\mupar_6$  & Rear end depression  & -15 mm & 0 mm \\
\hline
\hline
\end{tabular}
\end{table}

\begin{figure}[ht!]
  \centering
\begin{tabular}{|m{0.305\textwidth} | m{0.305\textwidth} | m{0.305\textwidth}|}
\hline
   \includegraphics[width=0.31\textwidth, clip, trim={0cm 0cm 0cm
   -2cm}]{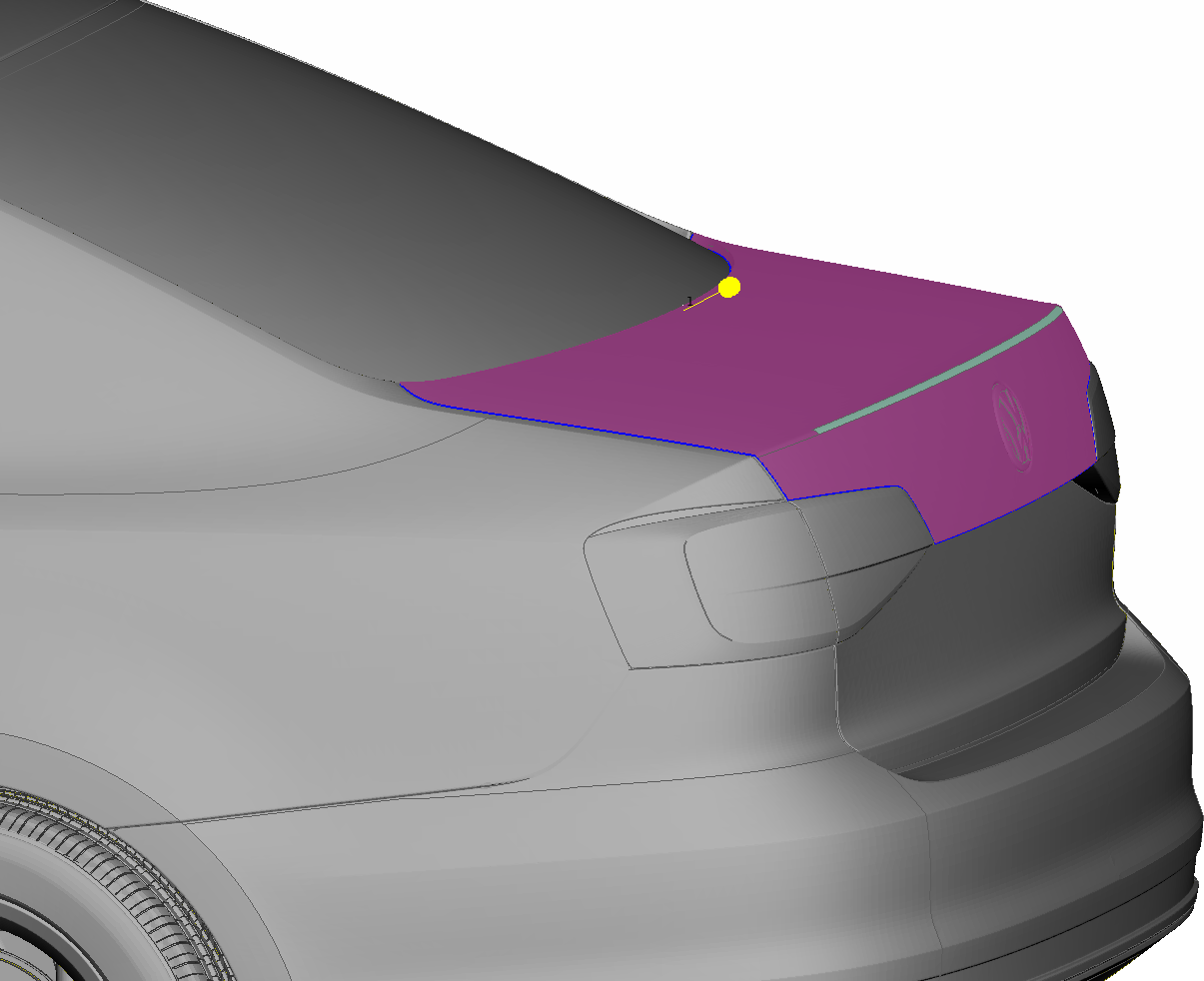} &
   \includegraphics[width=0.31\textwidth, clip, trim={0cm 0cm 0cm
   -2cm}]{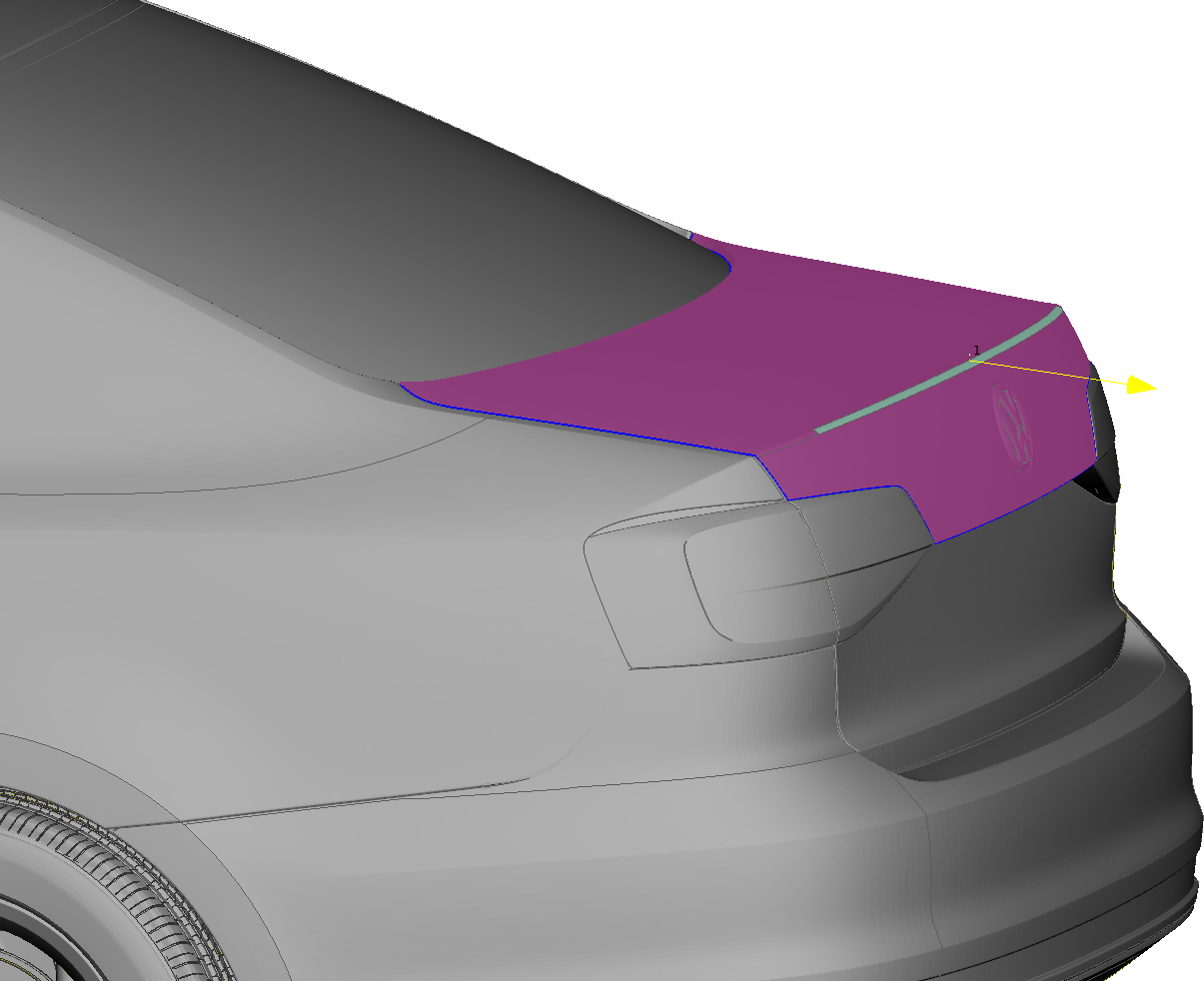} &
   \includegraphics[width=0.31\textwidth]{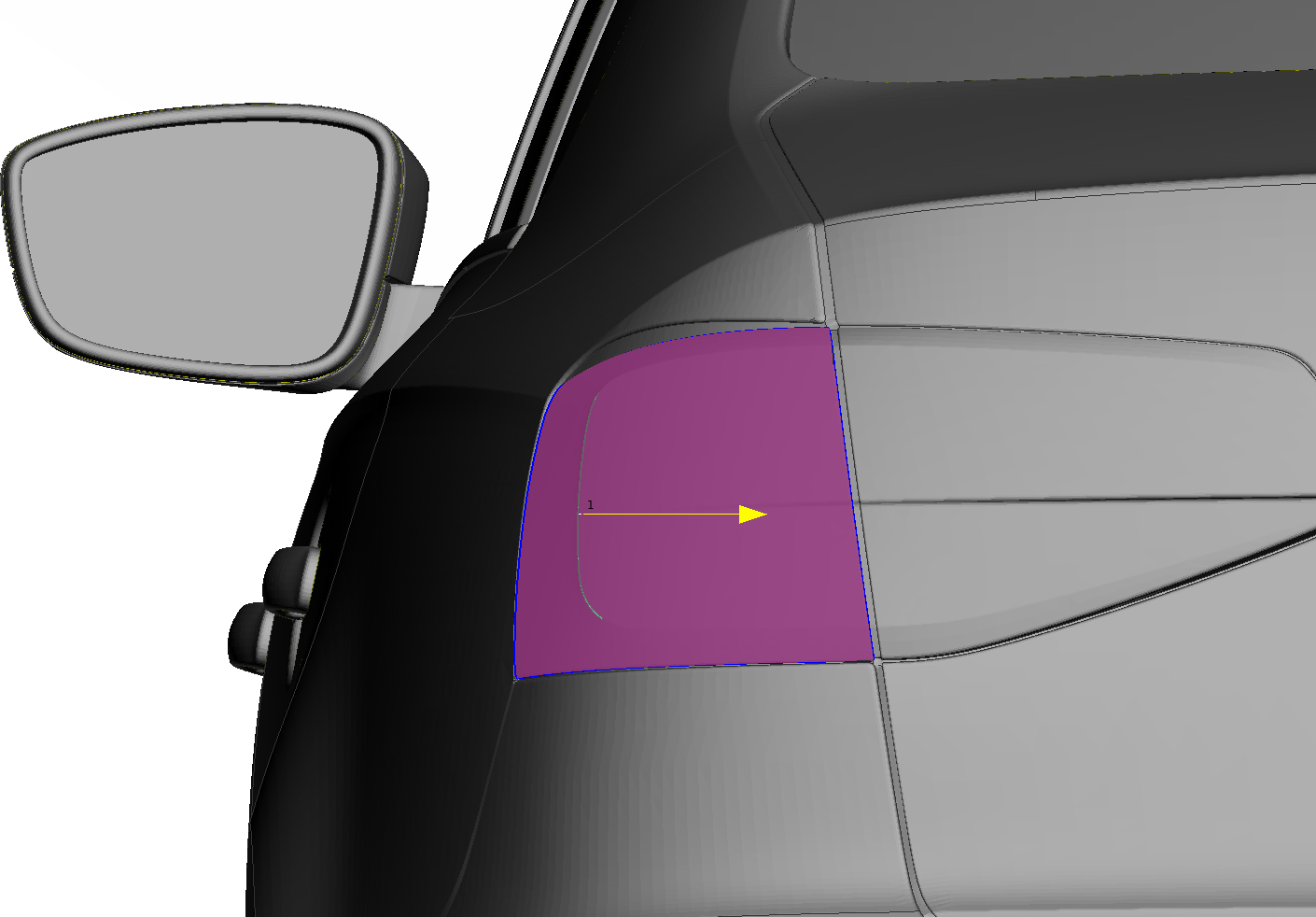}  \\
 \begin{center} Spoiler $Y$ Angle \end{center} & \begin{center} Spoiler Slide
 Translation \end{center} & \begin{center} Tail Light $Y$ Span \end{center}\\
 \hline
   \includegraphics[width=0.31\textwidth, clip, trim={0cm 0cm 0cm
   -2cm}]{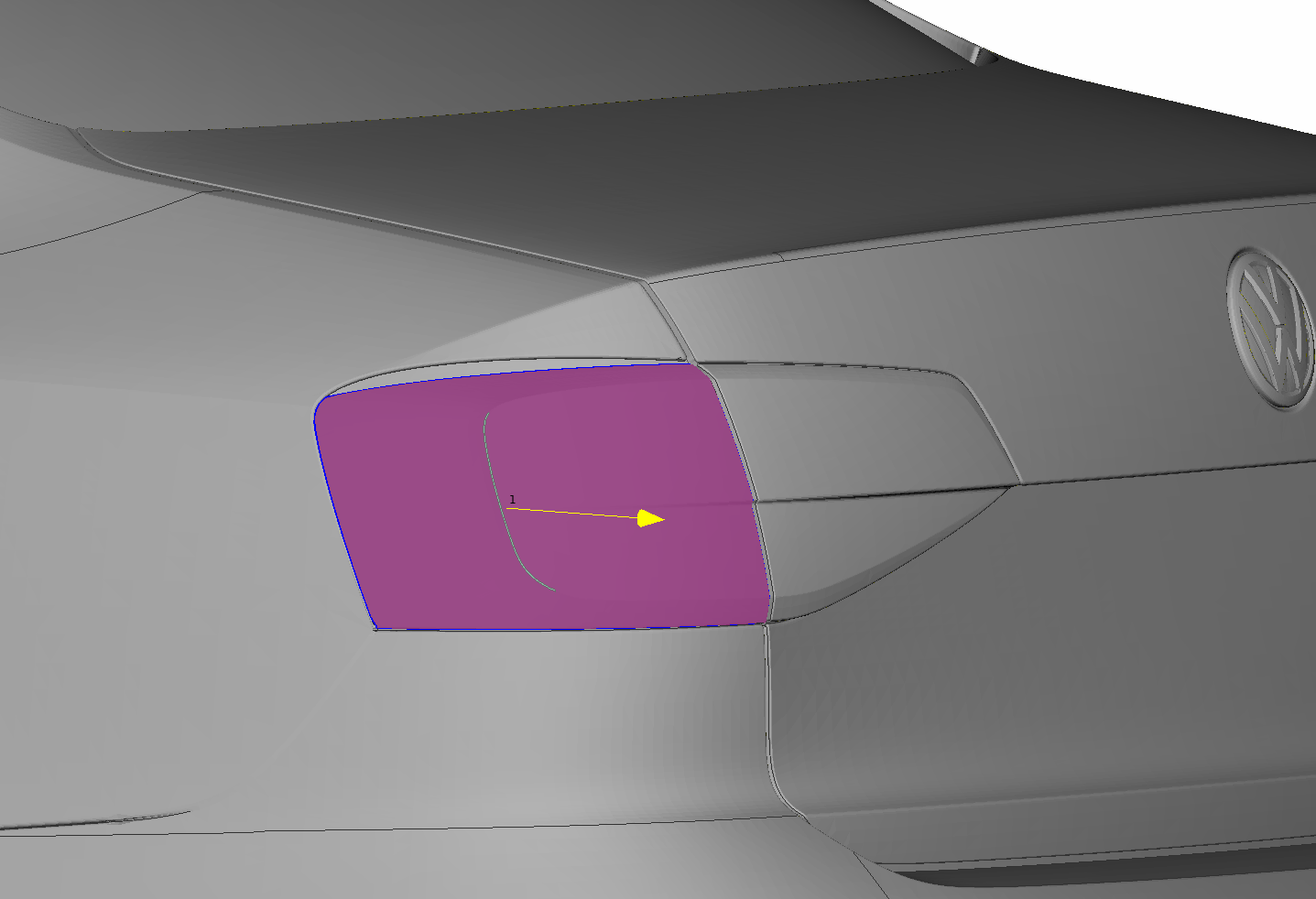} &
   \includegraphics[width=0.31\textwidth]{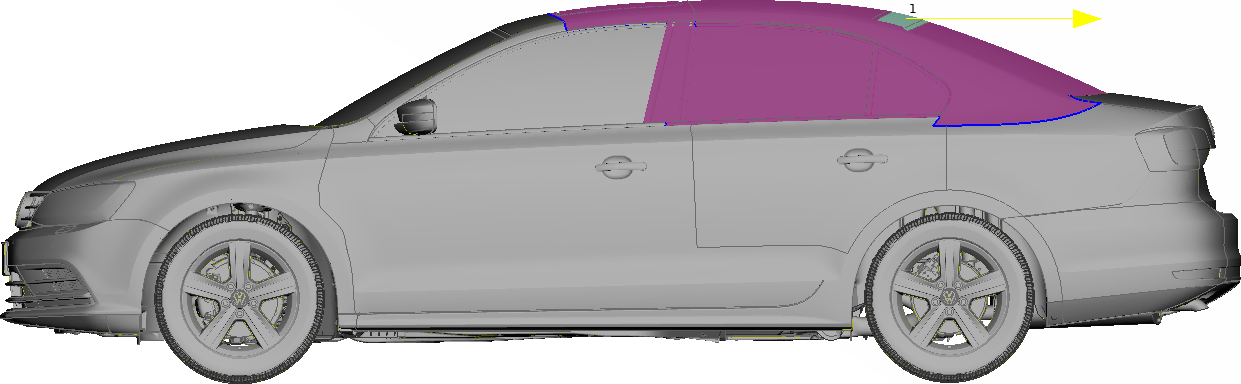} &
   \includegraphics[width=0.31\textwidth, clip, trim={0cm 0cm 0cm
   4.7cm}]{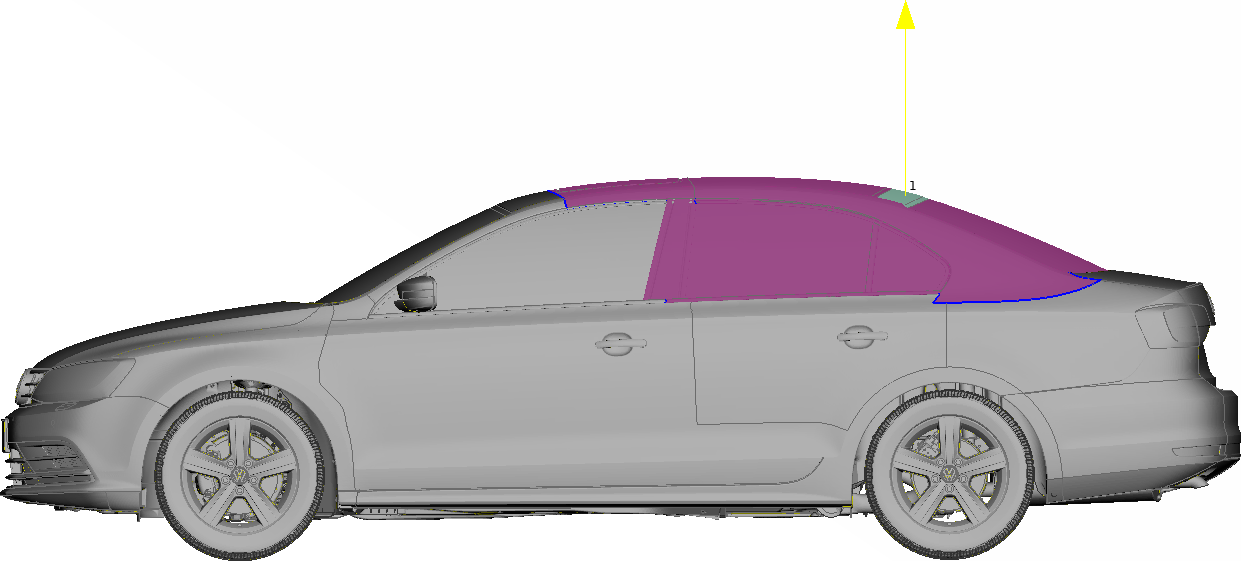}  \\
 \begin{center} Tail Light $X$ Translation \end{center} &  \begin{center} Rear
 Window $X$ Translation \end{center} & \begin{center} Rear Window $Z$
 Translation \end{center} \\
 \hline
   \includegraphics[width=0.31\textwidth, clip, trim={0cm 0cm 0cm
   1cm}]{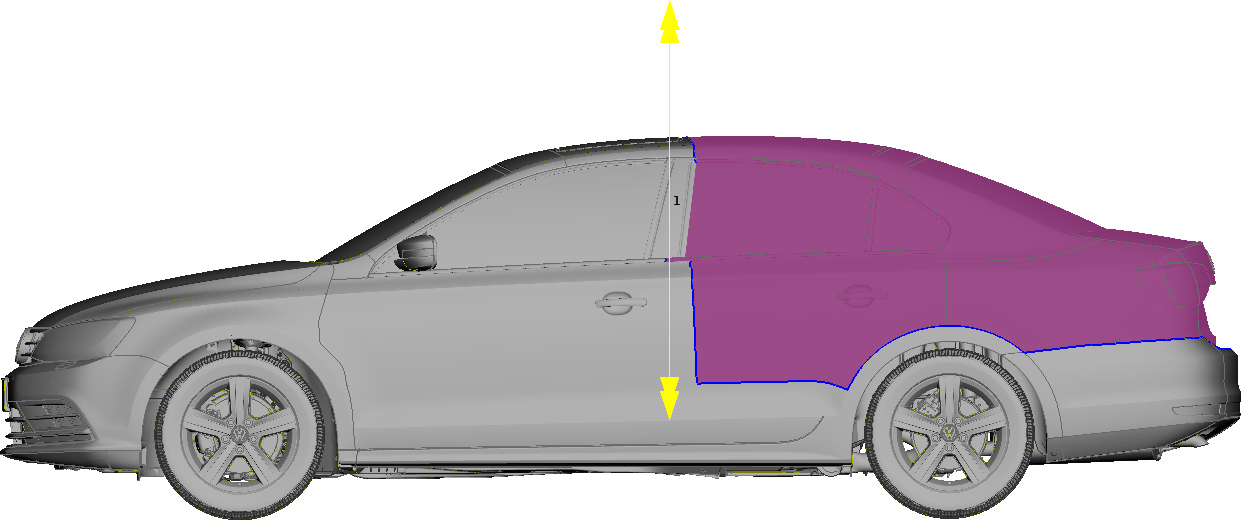} &
   \includegraphics[width=0.31\textwidth, clip, trim={0cm 0cm 0cm
   -2.7cm}]{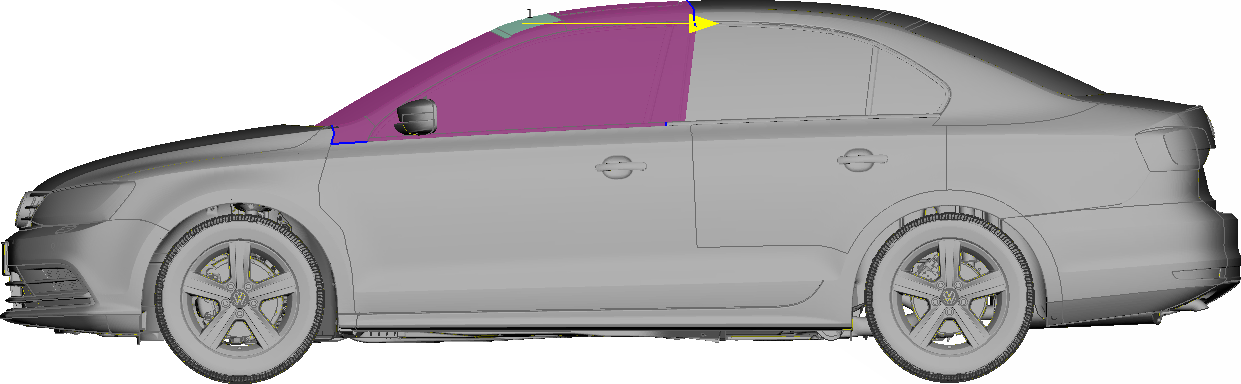} &
   \includegraphics[width=0.31\textwidth, clip, trim={0cm 0cm 0cm
   2cm}]{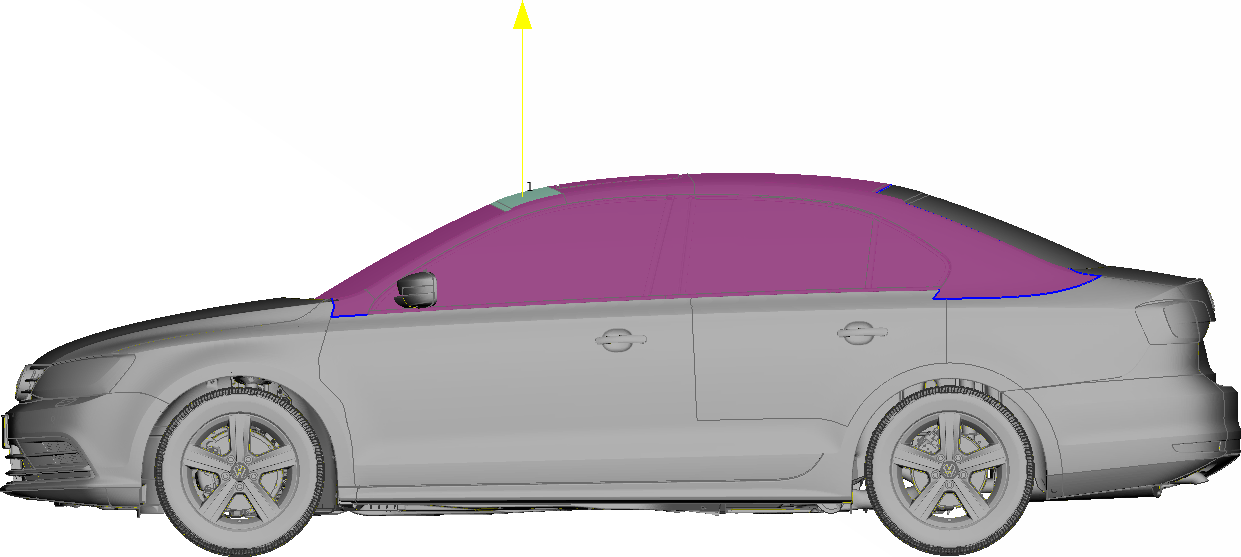}  \\
  \begin{center} Rear End Taper Ratio \end{center} & \begin{center} Front Window
  $X$ Translation \end{center} & \begin{center} Front Window $Z$ Translation
  \end{center} \\
  \hline
    \includegraphics[width=0.31\textwidth]{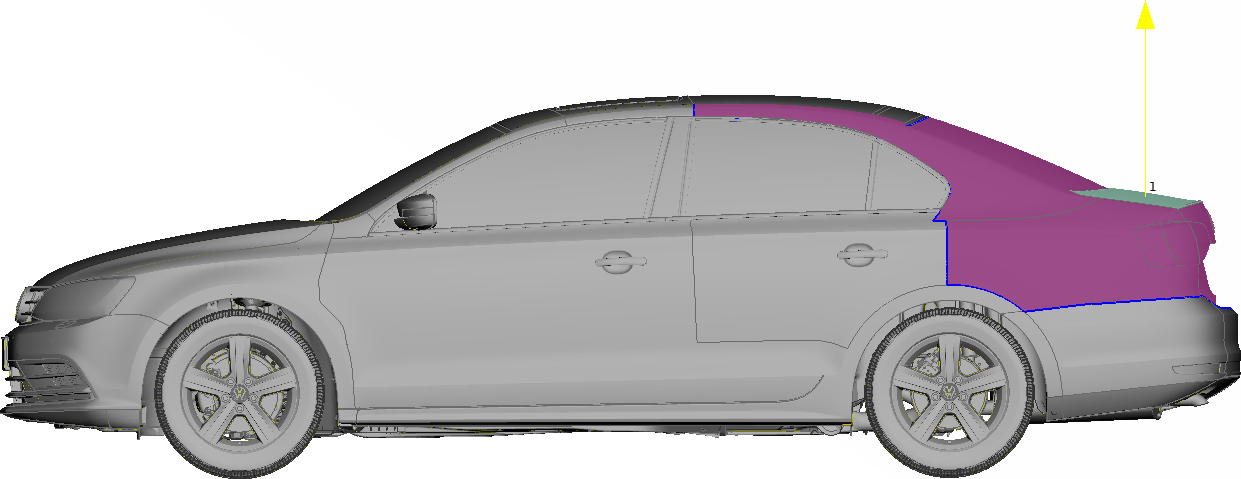} &
    \begin{center}\includegraphics[width=0.27\textwidth, clip, trim={0cm 0cm 0cm
    -2cm}]{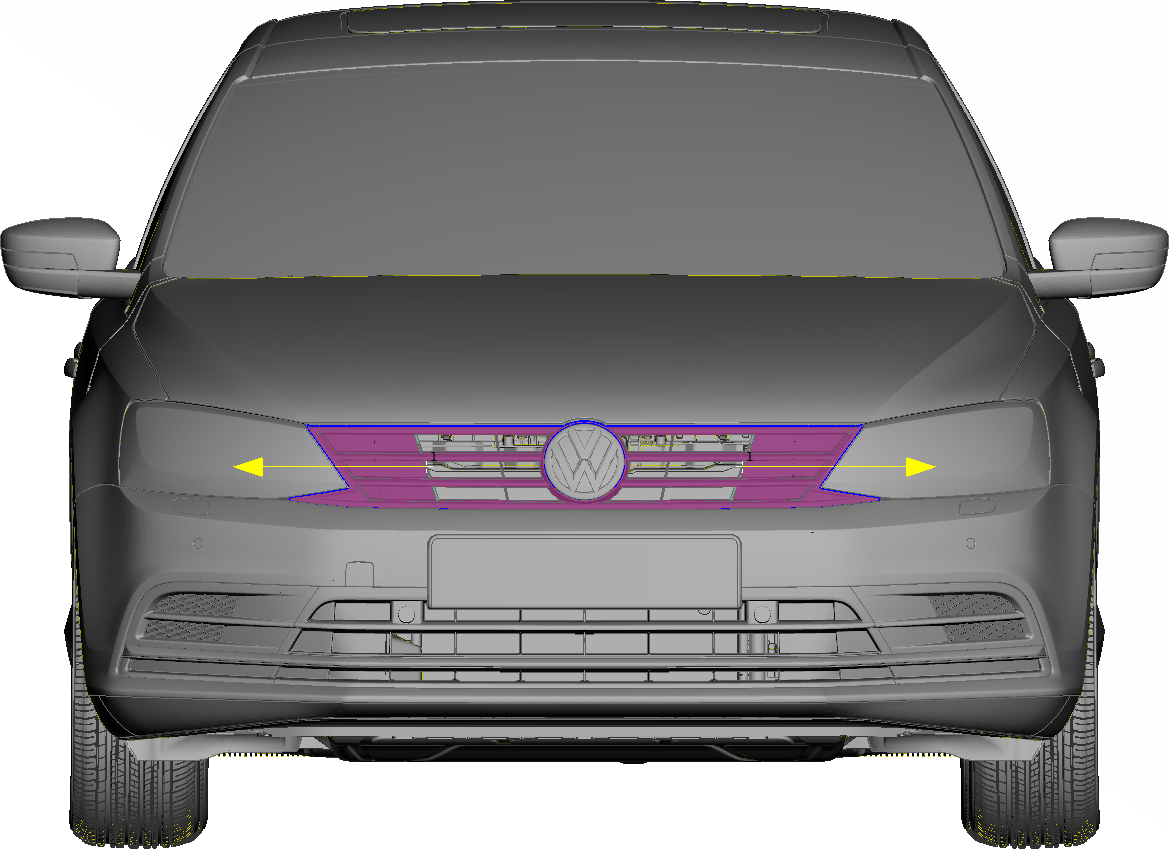} \end{center} &
    \includegraphics[width=0.31\textwidth, clip, trim={0cm 0cm 0cm
    -2cm}]{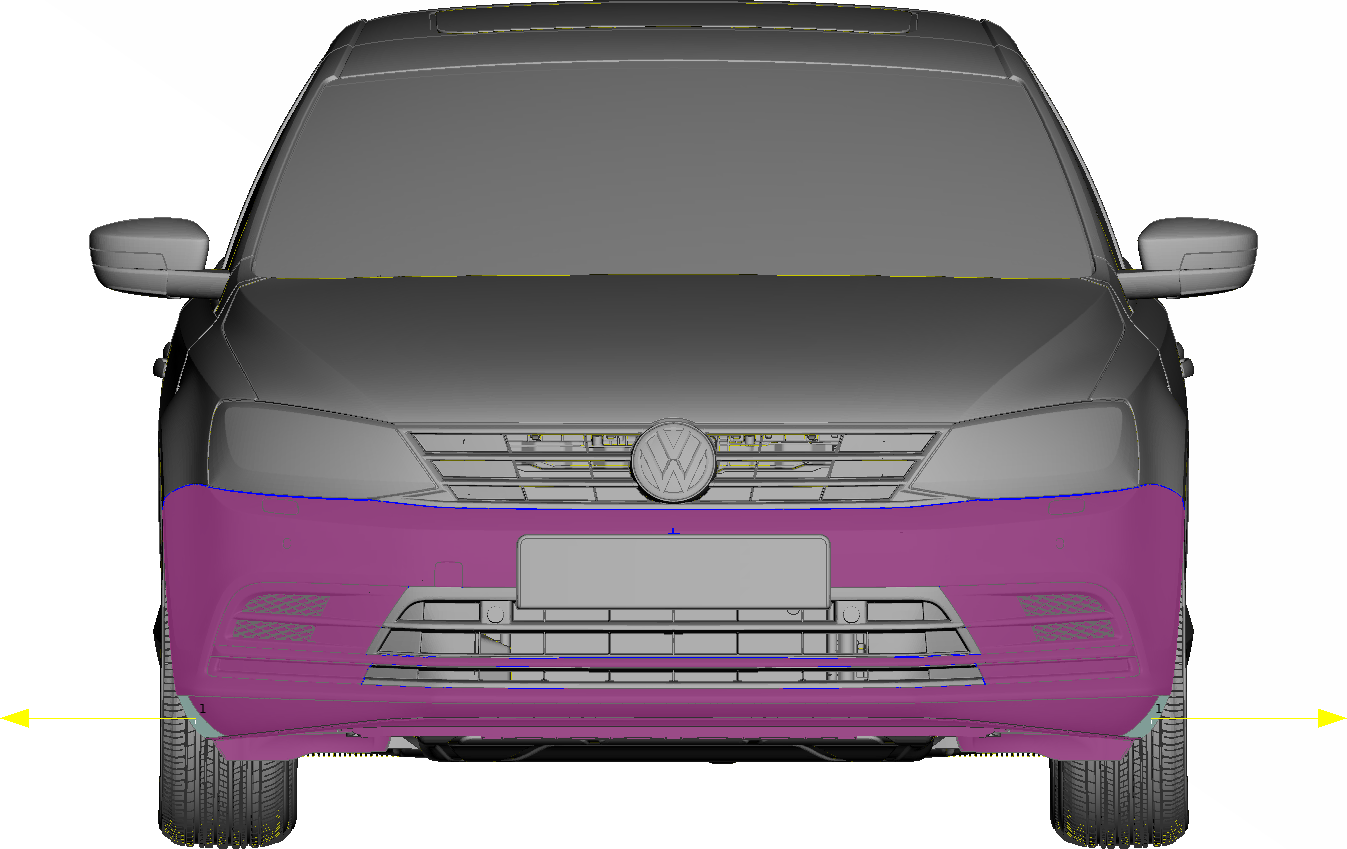}  \\
   \begin{center} Rear End $Z$ Translation \end{center} & \begin{center} Grill
   Slide Translation \end{center} & \begin{center} Bumper $Y$ Translation
   \end{center} \\
   \hline

\end{tabular}
  \caption{Affected areas by the geometrical parameters for the Jetta-12 test
    case. The ranges of each parameter can be gleaned from
    Table~\ref{tab:jetta_12_params}.}
  \label{fig:jetta_12_params}
\end{figure}

The second automotive test case (named hereafter {\bf Jetta-12}) is based on the
same car model and was created within the EC project UPSCALE~\cite{upscale}. The
parameterization consists of $12$ geometric modifications all around the vehicle
(see Figure~\ref{fig:jetta_12_params} and Table~\ref{tab:jetta_12_params}).
Besides the baseline geometry, a Sobol sequence of $300$ additional samples was
created and computed with OpenFOAM. To reduce the required computational budget
to an affordable amount, Reynolds-Averaged-Navier-Stokes (RANS) computations
were carried out instead of DDES runs. This allowed to use coarser meshes of
$52$M cells and resulted in $1700$ CPU-core-hours for a single run for the
$4000$ iterations, of which the last $1000$ were averaged to obtain the drag
coefficient $c_D$.

\begin{table}[htb]
\centering
\caption{Parameters' description of the Jetta-12 test
case.\label{tab:jetta_12_params}}
\begin{tabular}{ c l r r }
\hline
\hline
Parameter & Description & Lower bound & Upper bound \\
\hline
\hline
\rowcolor{Gray}
$\mupar_1$  & Spoiler $Y$ Angle  & -5.0° & 0.0° \\
$\mupar_2$  & Spoiler Slide Translation  & 0.0 mm & 30.0 mm \\
\rowcolor{Gray}
$\mupar_3$  & Tail Light $Y$ Span  & -15.0 mm & 5.0 mm \\
$\mupar_4$  & Tail Light $X$ Translation  & -10.0 mm & 10.0 mm \\
\rowcolor{Gray}
$\mupar_5$  & Rear Window $X$ Translation  & -100.0 mm & 100.0 mm \\
$\mupar_6$  & Rear Window $Z$ Translation  & -30.0 mm & 0.0 mm \\
\rowcolor{Gray}
$\mupar_7$  & Rear End Taper Ratio  & -1.0° & 3.0° \\
$\mupar_8$  & Front Window $X$ Translation  & -100.0 mm & 100.0 mm \\
\rowcolor{Gray}
$\mupar_9$  & Front Window $Z$ Translation  & -30.0 mm & 0.0 mm \\
$\mupar_{10}$  & Rear End $Z$ Translation  & -30.0 mm & 30.0 mm \\
\rowcolor{Gray}
$\mupar_{11}$  & Grill Slide Translation  & -50.0 mm & 50.0 mm \\
$\mupar_{12}$  & Bumper $Y$ Translation  & -20.0 mm & 20.0 mm \\
\hline
\hline
\end{tabular}
\end{table}

\begin{figure}[ht!]
    \centering
    \includegraphics[width=.85\textwidth]{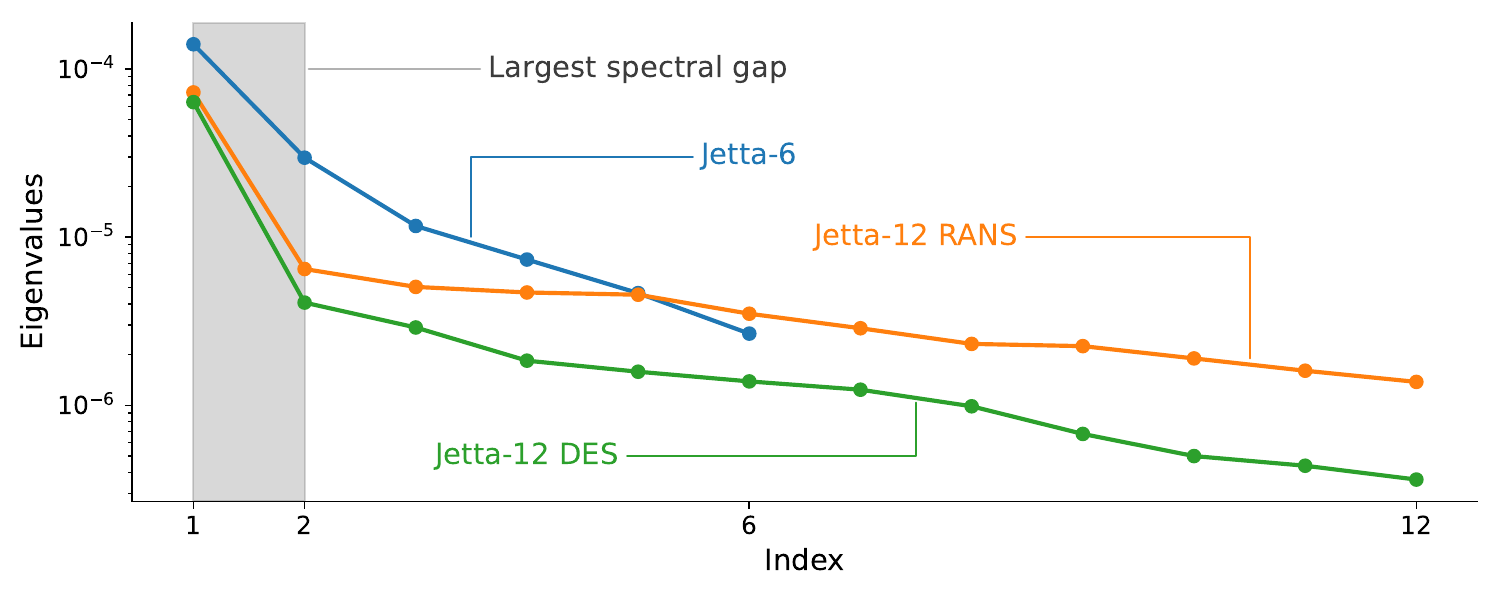}
    \caption{\RB{Eigenvalues decay of the covariance matrix of the
        gradients for the Jetta-6 and Jetta-12 test cases.}}
    \label{fig:eigs_jetta}
  \end{figure}

\RB{In Figure~\ref{fig:eigs_jetta} we depicted the eigenvalues decay for
  the automotive test cases. The largest spectral gap is always
  between the first and the second eigenvalue. This justifies the choice
  of a low-fidelity model built from a one-dimensional regression.}
  
\subsubsection{Multi-fidelity response surface design Jetta-6}
\label{subsubsec:Jetta-6}

In this test case, the low-fidelity model chosen is the response surface trained
on the active latent variables obtained with the NLL method: instead of
prolonging along the orthogonal directions the one-dimensional regression built
on the active subspace, a GPR is trained on the deformed high-fidelity inputs
$\{g_{\text{NLL}}(\mathbf{x}^{H}_i)\}_{i=1}^{N_{H}}\subset\Tilde{\mathcal{X}}$.
We remark that the map $g_{\text{NLL}}$ does not preserve in general convexity
of the domain $\mathcal{X}$ or orthogonality of the boundaries. Nonetheless,
this is not problematic for this application since we are not interested in
backmapping the active latent variables from $\Tilde{\mathcal{X}}$ to
$\mathcal{X}$, but only in forwarding the inputs from $\mathcal{X}$ to
$\Tilde{\mathcal{X}}$ and than evaluating the predictions with the GPR.

\begin{figure}[ht!]
  \centering
  \includegraphics[width=.85\textwidth]{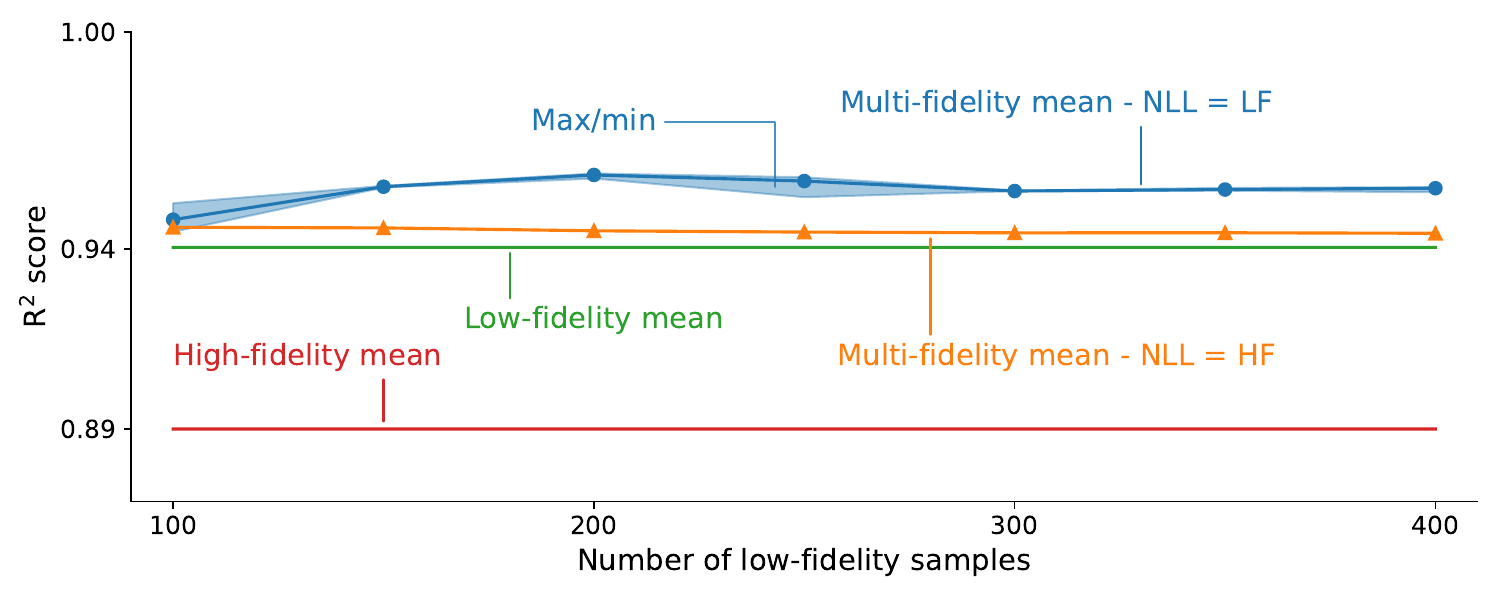}
  \caption{\RB{$R^2$ score evaluated on the $25$ test samples obtained from LHS on
  the domain $\mathcal{X}$, varying the number of LF samples. The mean $R^2$
  score over $10$ restarts of the training of the GPR is shown. For the MF also
  the minimum and maximum values are reported. The orange line
  identifies the results obtained by reversing the fidelities order,
  so the number of LF samples corresponds to the HF GPR. The LF and HF $R^2$ scores
  are not influenced by the number of additional LF samples.}}
  \label{fig:r2_scores_jetta6}
\end{figure}

The employed RevNet has $10$ layers. It was trained for $20000$ epochs on a
dataset of $76$ training samples and $25$ test samples, with ADAM stochastic
optimization method~\cite{kingma2017adam}, with an initial learning rate of
$0.03$. The high-fidelity samples were obtained with LHS method. The
architecture is implemented in PyTorch~\cite{NEURIPS2019_9015} inside the
ATHENA~\cite{romor2020athena} Python package. We perform a study on the number of additional LF samples, distributed uniformly
on the domain, from $100$ to $400$ with a step of $50$. The results are shown in
Figure~\ref{fig:r2_scores_jetta6}.

The maximization of the log-likelihood is performed with $10$ restarts for the
HF and LF models, and $100$ restarts for the MF model, all inside GPy
optimization algorithm. All training procedures are moreover restarted $10$
times, testing the stability of the optimization process for each fidelity
model. This is done in order to show, in
Figure~\ref{fig:r2_scores_jetta6} with blue lines,
that the MF training presents some small instabilities with respect to the HF
and LF training, as expected. The LF and HF models are designed over the same HF
inputs-outputs datasets, so they are not influenced by the additional LF
samples.

\begin{remark}[Reversing the fidelities order]
  \label{rmk:reversing fidelities order}
  A natural question that may arise regards the correct ordering of the HF and
  LF models in the MF when the accuracy is higher for the LF as in
  Figure~\ref{fig:r2_scores_jetta6}. We perform a study with respect to the
  number of additional samples from the HF GPR (not from the numerical
  simulations), now the lowest fidelity in the MF model. Moreover, in order to
  reach a desirable accuracy, we add
  to each of the 2 levels of fidelity of the MF model $200$ uniformly sampled
  input-output pairs: the highest fidelity is the NLL GPR built with $76 + 200$
  training data; the lowest fidelity is the HF GPR built with training data
  equal to $76$ from numerical simulations $+ 200$ fictitiously from the HR GPR
  (not from numerical simulations) $+$ additional samples from $100$ to $400$
  with a step of $50$ from HR GPR (not numerical simulations). The results are
  reported in Figure~\ref{fig:r2_scores_jetta6} with orange lines (NLL
  = HF). The $R^2$ score is
  lower than the previous case. Generally, the ordering of the fidelities
  depends on the availability of data and the cost for obtaining them.

\end{remark}

\begin{figure}[ht!]
  \centering
  \includegraphics[width=1\textwidth]{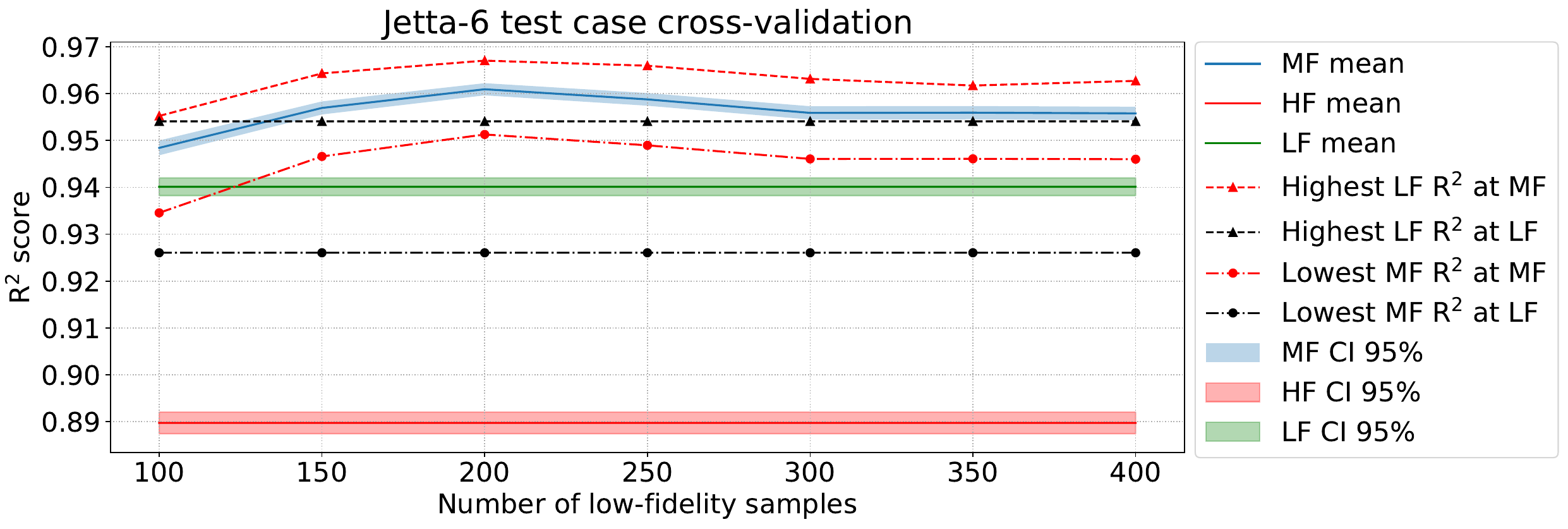}
  \caption{Cross-validation with leave-one-out strategy and confidence bounds at
  $95\%$. The labels \textbf{lowest MF $\mathbf{R^2}$ at LF} stands for the
  $R^2$ score of the batch associated to the lowest $R^2$ score for the MF model
  in the CV procedure, but evaluated at the predictions of the LF model. The
  other labels are analogue.}
  \label{fig:multifidelity_jetta6_cv}
\end{figure}

We also perform cross-validation (CV) with leave-one-out strategy for the
Jetta-6 test case to assess the robustness of the result with respect to the
test dataset, in Figure~\ref{fig:multifidelity_jetta6_cv}. We reported the mean
and confidence intervals at $95\%$ among the $25$ batches of the leave-one-out
strategy for a test set of $25$ samples: each batch has $24$ test samples. For
each abscissa, the batches corresponding to the lowest $R^2$ score for the MF
and highest $R^2$ score for the LF are found, so that with respect to these two
selected batches the $R^2$ scores of the LF and MF models, respectively, can be
computed and compared: we want to remark that batch-wise the MF $R^2$ score is
always higher to the LF $R^2$ score.

\subsubsection{Multi-fidelity response surface design Jetta-12}
\label{subsubsec:Jetta-12-RANS}

For this test case with additional $6$ parameters with respect to the previous
one, for a total of $12$, a one-dimensional NLL response surface does not
perform better than a one-dimensional AS response surface, so we preferred the
latter as LF model. In this case we also added Gaussian noise at each fidelity level in
order to achieve a better accuracy, loosing the Markov property, see
Remark~\ref{rmk:Markov}. Moreover, to avoid overfitting we restricted the variance of the Gaussian noise to the interval $[0.01, 0.1]$ at each fidelity of the multi-fidelity model.

\begin{figure}[ht!]
  \centering
  \includegraphics[width=1.\textwidth]{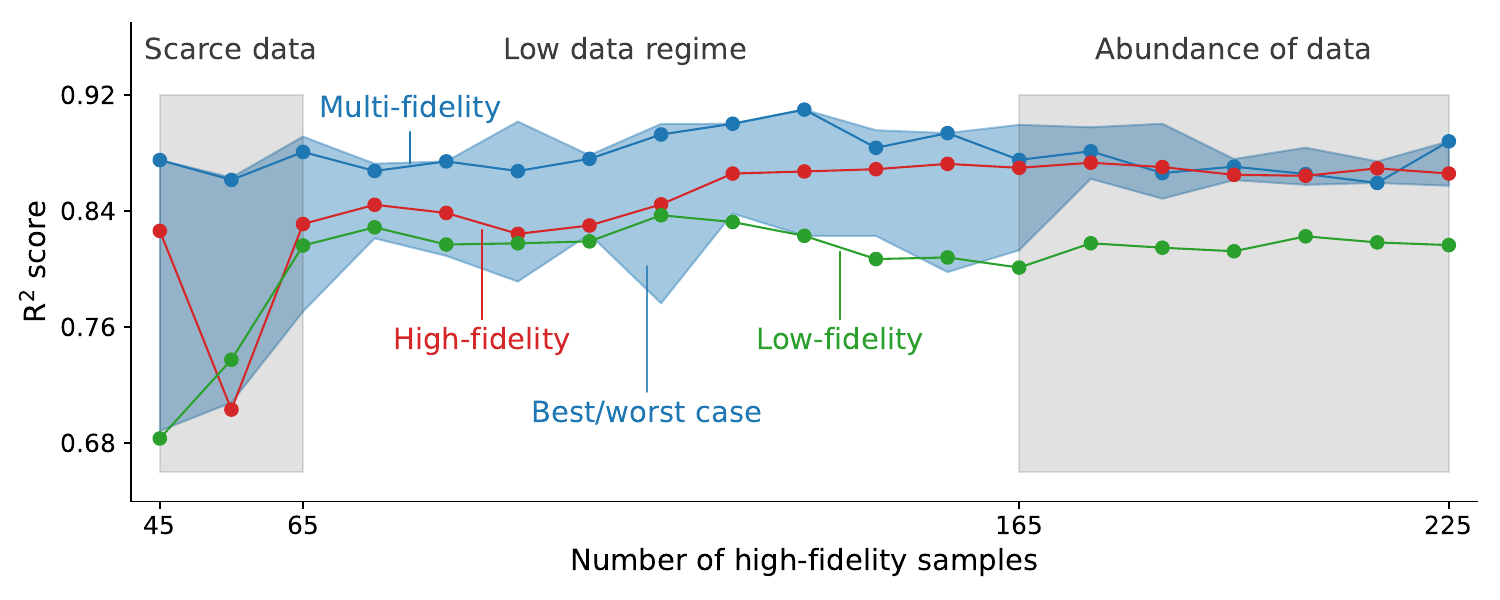}
  \caption{\RA{Jetta-12:} $R^2$ score evaluated on the $51$ test samples obtained from LHS on
  the domain $\mathcal{X}$, varying the number of HF samples. The mean $R^2$
  score over $10$ restarts of the training of the GPR is shown. For the MF model
  the minimum and maximum values are shown, differently from the HF and LF
  models, since the perturbations are not sensible. The markers
  associated to the MF and HF models represent the $R^2$ scores on the
  test set of the best HF and MF models with respect to an independent
  validation set of $25$ samples.}
  \label{fig:multifidelity_jetta12}
\end{figure}

We perform a study on the number of high-fidelity samples from $45$ to $225$,
obtained from a Sobol' sequence. The test set has $51$ samples obtained with LHS instead. The number of additional LF samples is $100$. The results are reported
in Figure~\ref{fig:multifidelity_jetta12}. As for the Jetta-6 test case, we perform $10$ outer training restarts for the LF, HF, and MF models: the $100$ additional LF samples are resampled every time. Moreover, the
optimization procedures of the GPRs are restarted $10$ times for the
LF, HF, and MF model. We employed also a validation dataset of
additional independent $25$ samples from the continuation of the
Sobol' sequence: the markers in Figure~\ref{fig:multifidelity_jetta12}
correspond to the best HF and MF models with respect to the validation
set. We also report maximum and minimum $R^2$ scores for the outer
loop training restarts of the MF model to show that the validation
process is fairly effective, at least when employing $45$ to $155$
high-fidelity samples. \RA{We emphasized in the plot the three
  distinct areas corresponding to the scarce data, low data, and
  abundance of data regimes.}

It can be seen a gain of around $4\%$ on average on the $R^2$ score of the MF model, with
respect to the other two, in the abscissae range from $45$ to $155$. This time the procedure is much less stable with
respect to the optimization process, probably due to the higher dimension of the
input space. The decreasing behaviour of the $R^2$ score of the MF models
from the abscissa $135$ to $225$ can be ascribed to the prevalence of the HF model: in this case the LF model influences less the predictions of the MF model, which are more stable and close to the HF ones. The low HF $R^2$ score at abscissa $55$ is almost constant for each outer training step and is not related to overfitting, but it can be associated to a high sensitivity of the regression when employing a small dataset relative to the problem at hand.

We want also to remark that in this test case the training of the LF, HF, and MF models takes less than $10$ minutes for each number of high-fidelity samples, with increasing costs from $45$ to $225$ samples, considering altogether the outer loop training restarts. Compared with the costs for a high-fidelity simulation,
the MF training cost is negligible.



\subsubsection{\RA{High-fidelity model choice for Jetta-12}}
\label{subsubsec:Jetta-12-RANS-DDES}

In principle, any model of the same phenomenon originated from a different
physical approximation, numerical method, or discretization, can be employed to
produce a multi-fidelity model. In the case of the Jetta-12 automotive testcase,
computations can be carried out with the more accurate DDES runs, as in the
Jetta-6 testcase. Then, we have 3 models at our disposal: the response surfaces built
on the DDES outputs, RANS outputs, or AS predictions.

We train the DDES-AS and DDES-RANS two-fidelity models as described in Section~\ref{sec:multifidelity}, and consider also the DES and AS
single fidelity models. We compute $75$
DDES training input-output pairs, and $50$ DDES test input-output pairs,
both sampled with LHS. Since the DDES simulations represent the highest
fidelity, the DDES test samples will be used to evaluate the $R^2$
scores of all the other single and multi-fidelity models considered. All the
$300$ RANS training data available from the previous test case will be employed for the DDES-RANS model, and
$100$ additional LHS sampled input-output pairs obtained from the AS response surface will
be used to train the DDES-AS two-fidelity model.

\begin{figure}[ht!]
  \centering
  \includegraphics[width=1\textwidth]{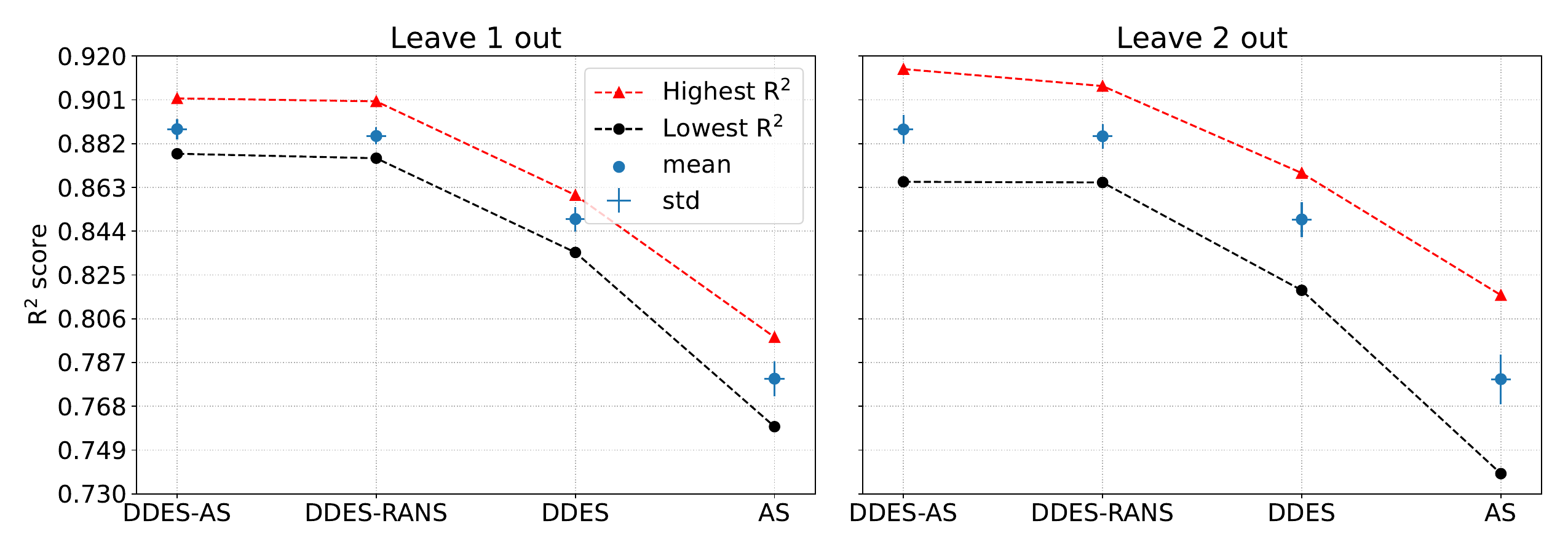}
  \caption{Cross-validation (CV) with leave-one-out (Left) and
  leave-two-out strategy (Right) on the test set. The labels
  \textbf{Lowest $\mathbf{R^2}$}
  and \textbf{Highest $\mathbf{R^2}$} stand
  for the $R^2$ score of the batch associated to the lowest and highest $R^2$
  score, respectively. The mean and standard deviation shown (std) are computed with
  respect to the sets of $50=\binom{50}{1}$ and $1225=\binom{50}{2}$
  cross-validation batches, respectively.}
  \label{fig:multifidelity_jetta12_3_fid_cv}
\end{figure}

\begin{figure}[ht!]
  \centering
  \includegraphics[width=0.85\textwidth]{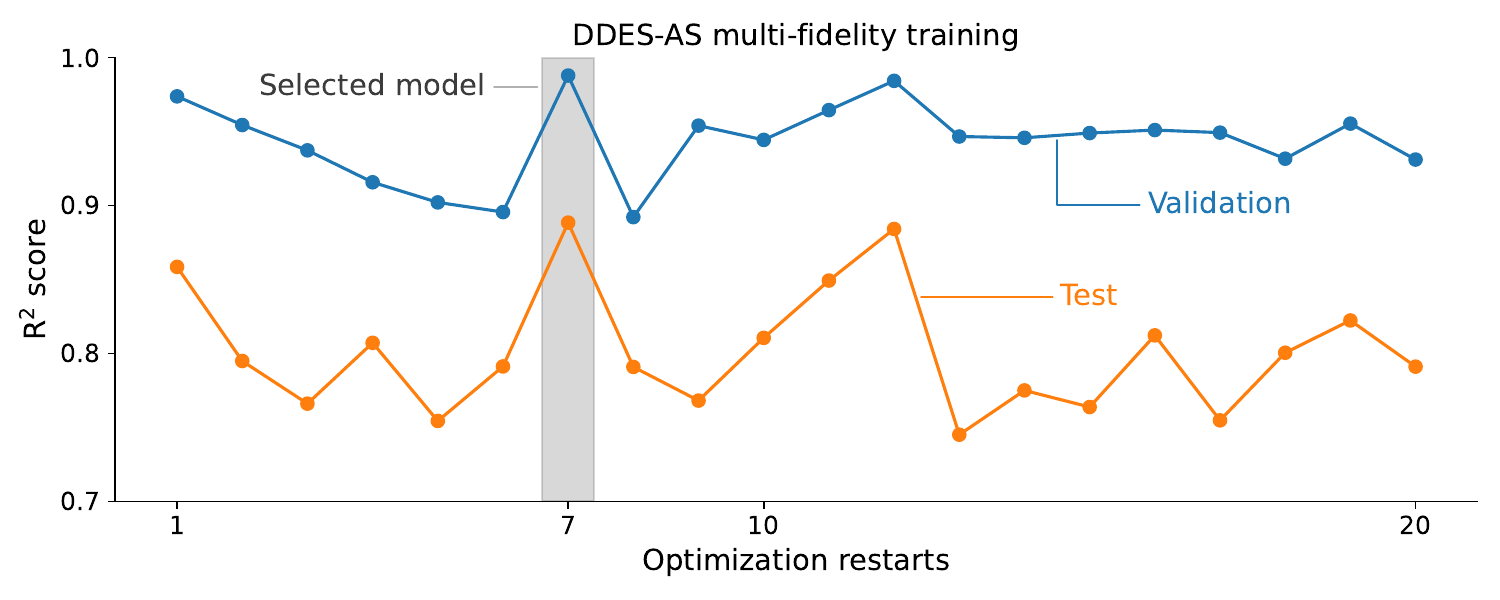}
  \caption{\RB{Validation process of the DDES-AS multi-fidelity model
      over $20$ outer training restarts changing every time the
      additional $100$ additional low-fidelity samples. The selected
      multi-fidelity model corresponds to abscissa $7$.}} 
  \label{fig:multifidelity_training}
\end{figure}

The $R^2$ errors on the test set are reported in
Figure~\ref{fig:multifidelity_jetta12_3_fid_cv}. Also in this case we use
cross-validation with leave-one-out and leave-two-out strategy to assess the robustness of the results with respect to the test
set: we show the mean, minimum, maximum, and standard deviation (std) with respect to the sets of
$50=\binom{50}{1}$ and $1225=\binom{50}{2}$ cross-validation batches. When the two models are integrated in the
DDES-AS MF model, the accuracy sensibly rises as observed in the previously. Only for the MF model, we extracted $10$ validation samples from the $75$ training dataset, so we trained it with exactly $65$ samples and selected the best model looking at the $R^2$ score of the validation set. The effectiveness of the validation procedure is shown in Figure~\ref{fig:multifidelity_training}. Also in this case, we constrained the Gaussian noise levels of the MF model to belong to the interval $[0.01, 0.1]$.

\begin{figure}[ht!]
  \centering
  \includegraphics[width=.8\textwidth]{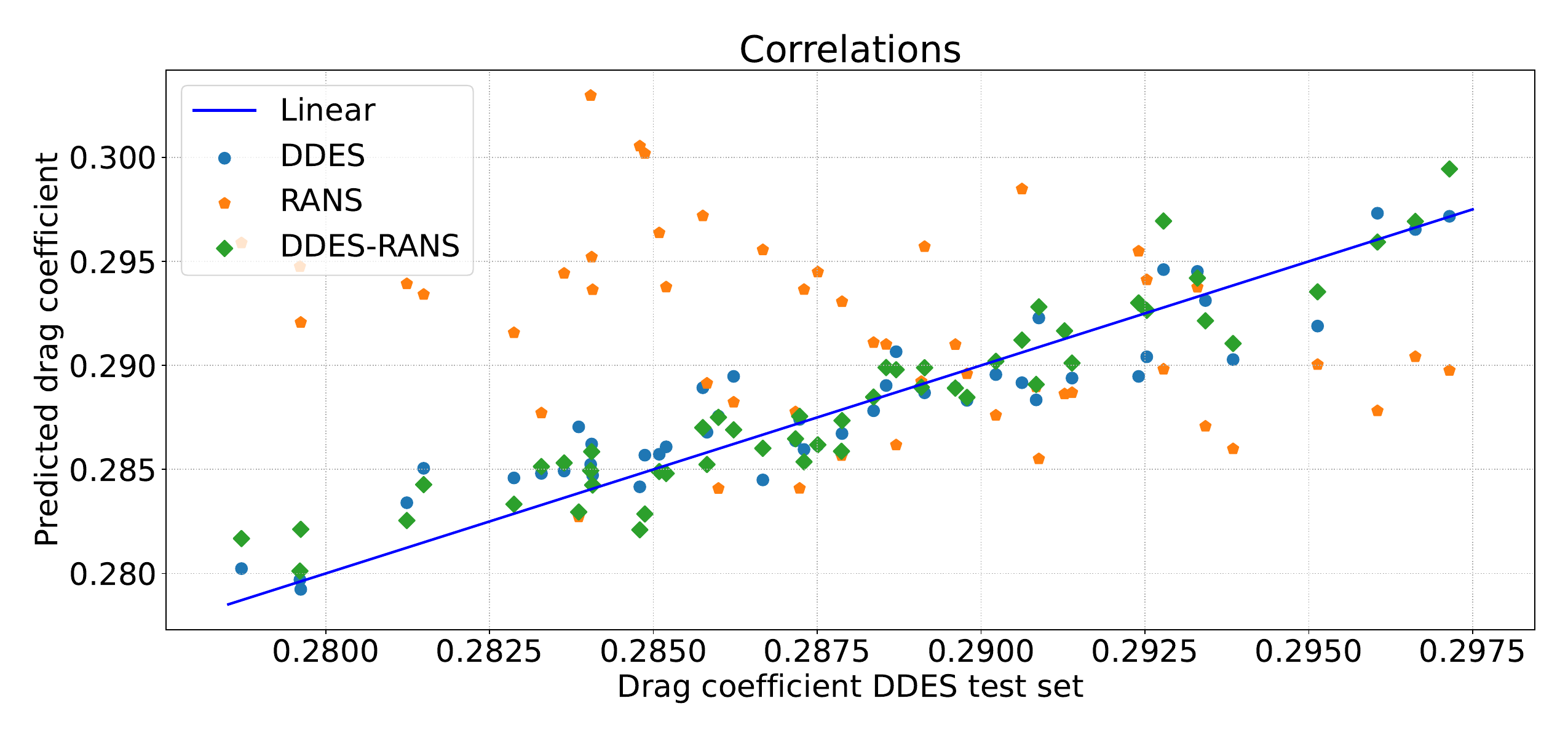}
  \caption{Comparison of the correlations between the predictions of
    the DDES model with the test DDES outputs, the correlations
    between the RANS model with the test DDES outputs, and the
    correlations between the DDES-RANS model with the test DDES
    outputs.}
  \label{fig:multifidelity_jetta12_3_fid_corr_rans}
\end{figure}

The accuracy is comparable to the DDES-RANS MF model, implying that the
AS response surface can indeed be used as a low-fidelity purely data-driven
model in the process of design of a multi-fidelity model, along with more standard models based on different physical or numerical
approximations of the phenomenon under study. It must be said that the RANS outputs are poorly correlated
with respect to the DDES as can be seen from
Figure~\ref{fig:multifidelity_jetta12_3_fid_corr_rans}: in fact the converged GPR built upon
the RANS training dataset have a mean $R^2$
score below $0$ on the DDES test set. Nonetheless, the multi-fidelity model
DDES-RANS achieves an accuracy higher than the single-fidelity DDES model.

\begin{figure}[ht!]
  \centering
  \includegraphics[width=.8\textwidth]{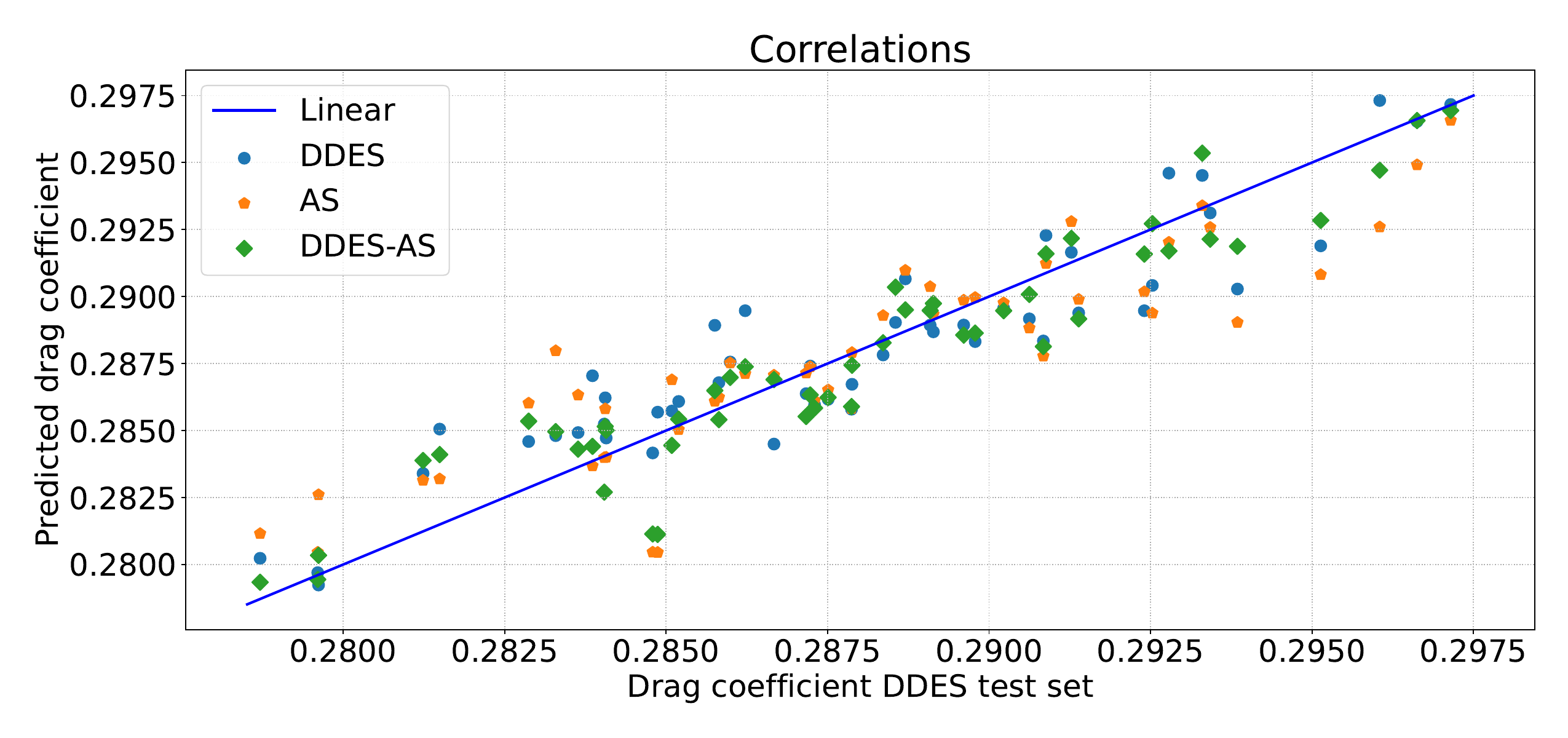}
  \caption{Comparison of the correlations between the predictions of the DDES model with the test DDES outputs and the correlations between the DDES-AS MF model with the test DDES outputs.}
  \label{fig:multifidelity_jetta12_3_fid_corr}
\end{figure}

\section{Conclusions and future perspectives}
\label{sec:conclusions}
The approximation of high-dimensional scalar quantities of interest is
a challenging problem in the context of data scarcity, which is
typical in engineering applications. We addressed this problems by
proposing a nonlinear multi-fidelity method which does not necessitate
the simulation of simplified models, but instead constructs a
low-fidelity surrogate introducing a low-intrinsic
dimensionality bias through active subspaces or nonlinear level-set
learning methods. Our approach is data-efficient since it extracts new
informations from the high-fidelity simulations. We construct different Gaussian
processes using the autoregressive scheme called NARGP. The proposed
multi-fidelity approach results in better approximation accuracy over
the entire parameter space as demonstrated with two benchmark problems
and an automotive application.

NARGP-AS was able to achieve better performance with respect to the
single-fidelity GP over the high-fidelity data, resulting in a
relative gain on the $R^2$ score around $3$--$5\%$  for the piston
model, and around $3$--$4\%$ for the Ebola model, depending on the
number samples used. NARPG-NLL was used for the Jetta-6 test case,
reaching an accuracy gain around $2\%$ with respect to the
low-fidelity model, and around $4\%$ with respect to the
high-fidelity model. We also presented a comparison switching the two
fidelities. Finally for the Jetta-12 test case we obtained a relative
gain on the $R^2$ score around $3\%$.

Future research lines should investigate the use of different active
subspaces-based methods, such as kernel AS~\cite{romor2022kas}, or local
AS~\cite{romor2021las}, which exploit kernel-based and localization
techniques, respectively. This multi-fidelity framework has also the
potential to be integrated with other reduced order modeling techniques~\cite{morhandbook2019,
chinestaenc2017,rozza2022book} to further increase the
accuracy in the resolution of parametric problems, especially for
high-dimensional surrogate-based optimization~\cite{tezzele2022multi}.

Mandatory for real applications is a model management strategy providing
theoretical guarantees and  establishing accuracy and/or convergence of
outer-loop applications. Some attempts towards multi-source Bayesian
optimization/Experimental design are being studied. Moreover increasing the
number of fidelities in the multi-fidelity model is a possible direction of
investigation, especially when the phenomenon of interest allows many cheap
low-fidelity approximations.

\section*{Acknowledgements}
This work was partially supported by the European Commission H2020 ARIA (Accurate ROMs for
Industrial Applications) project, by an industrial Ph.D. grant sponsored by
Fincantieri S.p.A. (IRONTH Project), by MIUR (Italian ministry for university
and research) through FARE-X-AROMA-CFD project, by the European Commission H2020
UPSCALE project (Upscaling product development simulation capabilities exploiting
artificial intelligence for electric vehicles, 824306), and partially funded by European
Union Funding for Research and Innovation --- Horizon 2020 Program --- in the
framework of European Research Council Executive Agency: H2020 ERC CoG 2015
AROMA-CFD project 681447 ``Advanced Reduced Order Methods with Applications in
Computational Fluid Dynamics'' P.I. Professor Gianluigi Rozza.

\bibliographystyle{abbrvurl}

\end{document}